\pdfoutput=1
\documentclass[a4paper,12pt,reqno]{amsart}
\usepackage[OT2,OT1]{fontenc} \newcommand\cyr
{
\renewcommand\rmdefault{wncyr} \renewcommand\sfdefault{wncyss} \renewcommand\encodingdefault{OT2} \normalfont
\selectfont
}
\DeclareTextFontCommand{\textcyr}{\cyr} \def\cprime{\char"7E }   

   \makeatletter
   \def\@settitle
{
\begin{center} \baselineskip14\p@\relax \LARGE\@title \end{center}
}
   \makeatother

\usepackage{graphics,color}
\usepackage{pstricks,epstopdf}
\usepackage{amsmath,amssymb,amsthm,graphicx,verbatim,psfrag,mathdots}
\usepackage[colorlinks=true]{hyperref}

\usepackage{appendix}
\usepackage[a4paper,margin=1.9cm,top=2.5cm,bottom=2.5cm,centering,vcentering]{geometry}

\graphicspath{{figs/}}
\theoremstyle{plain}
\newtheorem{theo}{Theorem}[section]
\newtheorem{lem}[theo]{Lemma}

\newtheorem{prop}[theo]{Proposition}

\theoremstyle{definition}

\newtheorem{rem}{Remark}

\newcommand{\pf}{\operatorname{Pf}}

\newcommand{\pff}{\operatorname{Pf}}

\newcommand{\one}{\raisebox{-1.5mm}{\scalebox{0.2}{\includegraphics{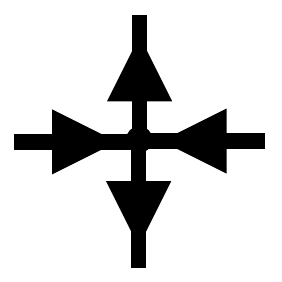}}}}

\newcommand{\nw}{\raisebox{-1.5mm}{\scalebox{0.2}{\includegraphics{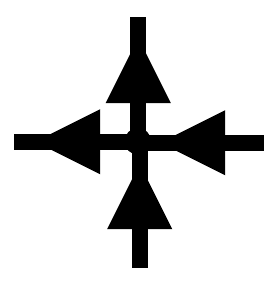}}}}
\renewcommand{\ne}{\raisebox{-1.5mm}{\scalebox{0.2}{\includegraphics{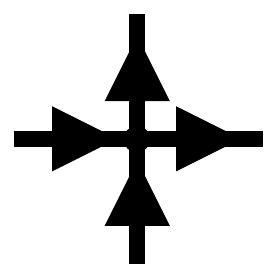}}}}

\newcommand{\qh}{\textup{QH}}

\newcommand\mycom[2]{\genfrac{}{}{0pt}{}{#1}{#2}}

\DeclareMathOperator{\av}{A_{V}}
\DeclareMathOperator{\avh}{A_{VH}}
\DeclareMathOperator{\avhpr}{A_{VHP}^R}
\DeclareMathOperator{\avhpc}{A_{VHP}^C}
\DeclareMathOperator{\ao}{A_{O}}
\DeclareMathOperator{\aoo}{A_{OO}}
\DeclareMathOperator{\aqt}{A_{QT}}
\DeclareMathOperator{\aqqt}{A_{qQT}}
\DeclareMathOperator{\aht}{A_{HT}}
\DeclareMathOperator{\an}{A}
\DeclareMathOperator{\avc}{A_{VC}}
\DeclareMathOperator{\avos}{A_{VOS}}
\DeclareMathOperator{\evh}{B_{VH}}
\DeclareMathOperator{\evhp}{B_{VHP}}
\DeclareMathOperator{\evos}{B_{VOS}}
\DeclareMathOperator{\sgn}{sgn}

\newcommand{\rightr}{\raisebox{-0.2mm}{\scalebox{1.0}{\includegraphics{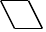}}}}
\newcommand{\vertr}{\raisebox{-0.2mm}{\scalebox{0.5}{\includegraphics{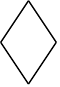}}}}

\numberwithin{equation}{section}

\usepackage[normalem]{ulem}

\begin{document}

\title[Refined Enumeration of ASMs]{Refined Enumeration of Symmetry Classes of Alternating Sign Matrices}
\author{Ilse Fischer and Manjil P. Saikia}
\address{Universit\"at Wien, Fakult\"at f\"ur Mathematik, Oskar-Morgenstern-Platz 1, 1090 Wien, Austria}
\email{\{ilse.fischer,~manjil.saikia\}@univie.ac.at}
\thanks{The authors acknowledge support from the Austrian Science Foundation FWF, START grant Y463 and SFB grant F50.}
\keywords{alternating sign matrices, six-vertex model, symmetry classes of alternating sign matrices, lozenge tilings of hexagons, non-intersecting lattice paths, symplectic group characters.}
\date{\today}

\newcommand\nnfootnote[1]{%
  \begin{NoHyper}
  \renewcommand\thefootnote{}\footnote{#1}%
  \addtocounter{footnote}{-1}%
  \end{NoHyper}
}
\nnfootnote{\emph{2010 Mathematics Subject Classification.} 05A05, 05A15, 05A19, 05E05, 15B35, 52C20, 82B20, 82B23.}  
\renewcommand{\thefootnote}{\fnsymbol{footnote}} 
\renewcommand{\thefootnote}{\arabic{footnote}}

\begin{abstract}
We prove refined enumeration results on several symmetry classes as well as related classes of alternating sign matrices with respect to classical boundary statistics, using the six-vertex model of statistical physics. More precisely, we study vertically symmetric, vertically and horizontally symmetric, vertically and horizontally perverse, off-diagonally and off-antidiagonally symmetric, vertically and off-diagonally symmetric, quarter turn symmetric as well as quasi quarter turn symmetric alternating sign matrices. Our results prove conjectures of Fischer, Duchon and Robbins.
\end{abstract}

\maketitle

\section{Introduction}\label{sec-one}

An \textit{alternating sign matrix} (ASM) of order $n$ is an $n \times n$ matrix with entries in the set $\{\pm 1,0\}$ such that all row and column sums are equal to $1$ and the non-zero entries alternate in each row and column. An example of an ASM of order $7$ is

\[
\begin{pmatrix}
0&0&0&1&0&0&0\\0&1&0&-1&0&1&0\\1&-1&0&1&0&-1&1\\0&0&1&-1&1&0&0\\0&1&-1&1&-1&1&0\\0&0&1&-1&1&0&0\\0&0&0&1&0&0&0
\end{pmatrix}.
\] These matrices, first introduced by David P. Robbins and Howard Rumsey, Jr. in the 1980s, have given rise to a lot of nice enumerative conjectures and results. Robbins, Rumsey and W. H. Mills \cite{asm-conj} conjectured that the number of ASMs of order $n$ is given by $$\prod_{i=0}^{n-1}\frac{(3i+1)!}{(n+i)!}.$$ This conjecture was later proved by Doron Zeilberger \cite{doron} and shortly after also by Greg Kuperberg \cite{asm-kuperberg}, using vastly different techniques. A very detailed description of this conjecture and Kuperberg's proof can be found in a book by David Bressoud \cite{proofs}. Another approach towards proving this result was developed by the first author \cite{ilse-1, ilse-2}.

 Richard Stanley \cite{stanley} suggested the study of symmetry classes of ASMs shortly after these objects were introduced. This led Robbins \cite{robbins-1, robbins-2} to conjecture various product formulas for the different symmetry classes of ASMs. The program of proving these formulas was accomplished by work of Kuperberg \cite{KuperbergRoof}, Soichi Okada \cite{OkadaCharacters}, A. V. Razumov and Yu. G. Stroganov \cite{RazumovStroganov,rs-odd-half,rs-odd-quarter}, and recent work of Roger E. Behrend,  Matja\v{z} Konvalinka and the first author \cite{bfk}.

It is easy to see that there is precisely one occurrence of $1$ in the first row of any ASM. This suggests the study of some refined enumerations of ASMs and symmetry classes thereof with respect to the position of this $1$. The refined enumeration of all ASMs with the position of the unique $1$ in the first row fixed was done by Zeilberger \cite{zeil}. In case of the symmetry classes of ASMs, this refinement has been studied for half-turn symmetric ASMs and quarter-turn symmetric ASMs. The refined enumeration for the first was provided by Razumov and Stroganov \cite{rs-odd-half}, while we prove a formula for the case of quarter-turn symmetric ASMs in this paper (in Section \ref{sec:qtsasm}) which was conjectured by Robbins \cite[Conjecture 3.6]{robbins-2}. In fact, we also prove a similar result for what are called quasi-quarter-turn symmetric ASMs that was conjectured by Philippe Duchon \cite{duchon}, but we defer this discussion until Sections \ref{sec:qtsasm} and \ref{sec:qqtsasm}.

With regard to vertically symmetric ASMs (VSASMs) and vertically and horizontally symmetric ASMs (VHSASMs), the refined enumeration with respect to the first row is trivial because the position of the $1$ is always fixed to be in the center. However, it turns out that in the second row of such matrices there are exactly two occurrences of $1$. The first author \cite{FischerOpFormulaVSASM} had conjectured a formula for the number of $(2n+1)\times (2n+1)$ VSASMs (the matrix above is an example of a VSASM), where the first 1 in the second row is in the $i$-th column, which we prove in this paper (in Section \ref{sec:vsasm}). No exact enumeration formula of such type has been conjectured for the case of VHSASMs. In Section \ref{sec:vhsasm}, we give generating function results for the refined enumeration of VHSASMs. We also prove refined enumeration results for other related matrices with vertical symmetry, namely vertically and horizontally perverse ASMs (VHPASMs) and vertically and off-diagonally symmetric ASMs (VOSASMs). Another class of ASMs with special symmetries are odd order off-diagonally and off-antidiagonally symmetric ASMs (OOSASMs). We prove generating function results for refined enumeration of such matrices as well.

By reflection, an ASM also has a unique $1$ in the last row as well as in the first and last column. Various mathematicians have studied related refinements, considering restrictions on a combination of two or more of the boundaries of an ASMs. These works include but are not limited to those by Arvind Ayyer and Dan Romik~\cite{AyyerRomik}, Behrend~\cite{behrend}, Fischer \cite{ilse-refined}, Fischer and  Romik~\cite{ilse-romik}, Razumov and Stroganov \cite{RazumovStroganov, rs-refined}, Romik and Matan Karklinsky~\cite{romik}, Stroganov~\cite{stroganov-izergin}, etc. using a variety of tools, but predominantly techniques that arise in statistical physics.

This paper is structured as follows: in Section \ref{asm} we shall briefly review the connection of ASMs with the six vertex model in statistical physics, which we will use for our proofs and then build our setup for the rest of the paper. In Section \ref{sec:vsasm} we study the six-vertex model and the associated partition function for VSASMs, and prove the first author's conjecture. In Section \ref{sec:vhsasm} we provide generating functions for the refined enumeration of VHSASMs, and in Section \ref{sec:vhp} we provide formulas for the refined enumeration of VHPASMs. Then, in Section \ref{sec:oos}, we give generating functions for the refined enumeration of OOSASMs, while in Section \ref{sec:vos} we give generating functions for the refined enumeration of VOSASMs. In Section \ref{sec:qtsasm} we prove formulas for the refined enumeration of quarter-turn symmetric ASMs which prove conjectures of Robbins \cite{robbins-2}, and, finally, in Section~\ref{sec:qqtsasm} we study the refined enumeration of quasi quarter-turn symmetric ASMs and prove Duchon's conjecture. We give some concluding remarks in Section \ref{sec:rem}. Some auxiliary results required in Sections \ref{sec:vhsasm}, \ref{sec:oos} and \ref{sec:vos} are proved in Appendix \ref{appen}.

\section{ASMs and the Six Vertex Model}\label{asm}

Kuperberg's proof of the ASM conjecture \cite{asm-kuperberg} is by exploiting a bijection between ASMs and a model in statistical physics, called the six-vertex model. In this section, we will explain this connection in brief.

\begin{figure}
\centering
\includegraphics{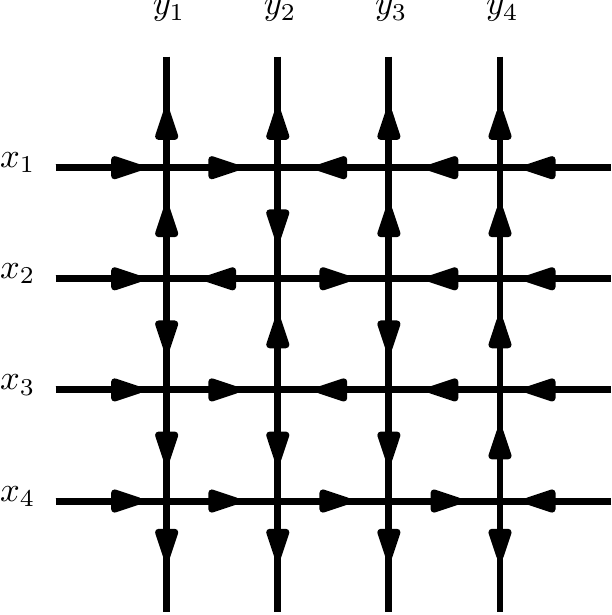}
\caption{Six Vertex Model with Domain Wall Boundary Condition.}
\label{fig:6-vertex}
\end{figure}

We consider a quadratic sub-region of the square grid such that the boundary vertices are of degree $1$, see Figure~\ref{fig:6-vertex}. A \textit{configuration} of a corresponding statistical physics model is an orientation on the edges of this graph such that both the in-degree and the out-degree of each degree $4$ vertex is $2$. If the orientations of the edges incident with vertices of degree $1$ are prescribed to be oriented inwards for horizontal edges and outwards for vertices edges (see Figure~\ref{fig:6-vertex}) then we say such a configuration satisfies the \textit{domain wall boundary condition}. Such a setup is in bijection with ASMs, and is called the \textit{six vertex model}, due to the six possibilities of assigning orientations to the edges incident with degree $4$ vertices to make the in-degree and out-degree equal to $2$. An example of such a configuration in given in Figure~\ref{fig:6-vertex}.

\begin{figure}
\begin{center}
\begin{tabular}{lcccccc}
ASM & $1$ & $-1$  & $0$  & $0$  & $0$  & $0$ \\ \hline \vspace{1mm}
six-vertex configuration & 
\raisebox{-0.6cm}{\scalebox{0.5}{\includegraphics{1.pdf}}} & 
\raisebox{-0.6cm}{\scalebox{0.5}{\includegraphics{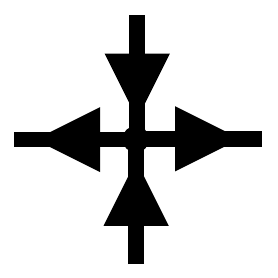}}} & 
\raisebox{-0.6cm}{\scalebox{0.5}{\includegraphics{ne.pdf}}} & 
\raisebox{-0.6cm}{\scalebox{0.5}{\includegraphics{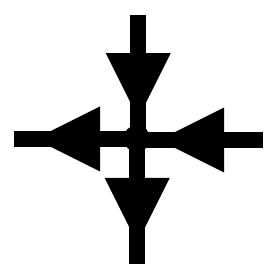}}} & 
\raisebox{-0.6cm}{\scalebox{0.5}{\includegraphics{nw.pdf}}} & 
\raisebox{-0.6cm}{\scalebox{0.5}{\includegraphics{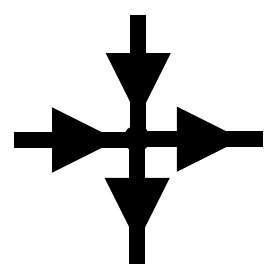}}} \\ \hline \vspace{1mm}
weight & $1$ & $1$ & $\frac{\sigma(q u)}{\sigma(q^2)}$ & $\frac{\sigma(q u)}{\sigma(q^2)}$ &
$\frac{\sigma( q \bar u)}{\sigma(q^2)}$ & $\frac{\sigma(q \bar u)}{\sigma(q^2)}$
\end{tabular}
\end{center}
\caption{\label{weights} Correspondence between ASMs and six-vertex configurations.}
\end{figure}

If we associate with each of the degree $4$ vertices in a six-vertex configuration a number as given in Figure \ref{weights}, then we obtain a matrix with entries in the set $\{\pm 1,0\}$. It is not difficult to see that such a matrix will be an ASM, and we actually get a bijection between ASMs and configurations of the six-vertex model. For example, the matrix associated with the configuration in Figure \ref{fig:6-vertex} is 
\[
\begin{pmatrix}
   0 & 1 & 0 & 0 \\
   1 & -1 & 1& 0 \\
   0 & 1 & 0 & 0 \\
   0 & 0 & 0 & 1
\end{pmatrix}.
\]

With such a model, there is an associated \textit{partition function} (say $Z$): it is the sum of the weights of all possible configurations (say $w( C )$) of the model. The \textit{weight of a configuration} $C$ is the product of all the weights of degree $4$ vertices (say $w_v$ for vertex $v$ in the graph), which will be defined now. We associate with each horizontal line in the grid a spectral parameter $x_i$ and with each vertical line a spectral parameter $y_j$ as marked in Figure \ref{fig:6-vertex}. The \emph{label} of a vertex which is intersected by lines associated with the parameters $x_i$ and $y_j$ is $x_i/y_j$. Each of the vertices will be assigned a weight as given in Figure \ref{weights}, where $u$ has to be replaced by the label and $q$ is a global parameter which we will specialize in the coming sections. We use the abbreviations 
$$\overline{x}=x^{-1}  \quad \text{and} \quad  \sigma(x)=x-\overline{x}$$ throughout this paper.

More formally, the above description says that the partition function is defined as  $$Z(n;\vec{x}; \vec{y})=\sum_{C \in {\mathcal C}}w( C )=\sum_{C \in {\mathcal C}}\prod_{\mycom{v\in C }{ \text{$v$ vertex of degree $4$}}}w_v,$$ where $\vec{x}=(x_1, x_2, \ldots, x_n)$, $\vec{y}=(y_1, y_2, \ldots, y_n)$ and $\mathcal C$ is the set of all possible configurations of the model we are considering. It is clear now that, if all the weights can be made equal to $1$ in all configurations of the six-vertex model by specializing $\vec{x}, \vec{y}$ as well as $q$, then the partition function will just count the number of ASMs of a given size. The choices of weights in Figure \ref{weights} are appropriate to make this happen. This was essentially the approach used by Kuperberg \cite{asm-kuperberg} to prove the ASM conjecture. 

We note that, rotating any of these vertices would result in changing the weights of the vertices in a way that the spectral parameters appearing in the weights would become their inverses. This will be used in the subsequent sections without commentary.

\section{Vertically Symmetric ASMs}\label{sec:vsasm}

ASMs that are invariant under the reflection in the vertical symmetry axis only exist for odd order: this follows because an ASM has a unique $1$ in the top row which then has to be situated in the central column of an ASM with vertical symmetry. This implies in particular that the top row of a VSASM is fixed. It is also not hard to see that the central column of a VSASM is always $(1,-1,1,-1, \ldots,1)^T$. The second row of a VSASM contains precisely two $1$'s, which are symmetrically arranged left and right of the central $-1$. Thus, the second row of a VSASMs is determined by the position of the first $1$ in the second row. The first author \cite{FischerOpFormulaVSASM} conjectured the following formula for the number of $(2n+1) \times (2n+1)$ VSASMs that have the first $1$ in the second row in position $i$:
$$
\frac{(2n+i-2)! (4n-i-1)!}{2^{n-1} (4n-2)! (i-1)! (2n-i)!} \prod_{j=1}^{n-1} \frac{(6j-2)! (2j-1)!}{(4j-1)! (4j-2)!}=: \av(2n+1,i).
$$
In this section, we prove this conjecture.

From the discussion in the previous paragraph it is clear that $(2n+1) \times (2n+1)$ VSASMs correspond to $2n \times (n+1)$ matrices with entries in $\{\pm 1,0\}$ that have the following properties.
\begin{enumerate}
\item The non-zero entries alternate in each row and column.
\item All column sums are $1$ except for the last column which is always $(1,-1,1,\ldots,-1)^T$.
\item The first non-zero entry of each row and column is $1$.
\end{enumerate}
The $6 \times 4$ matrix with these properties that corresponds to the VSASM from Section \ref{sec-one} is
$$
\begin{pmatrix}
0 & 0 & 0 & 1\\
0 & 1 & 0 & -1  \\
1 & -1 & 0 & 1  \\
0 & 0 & 1 & -1  \\
0 & 1 & -1 & 1 \\
0 & 0 & 1 & -1 
\end{pmatrix}.
$$
(We deleted the bottom row as well as the last $n$ columns.) 

\begin{figure}
\scalebox{0.8}{\includegraphics{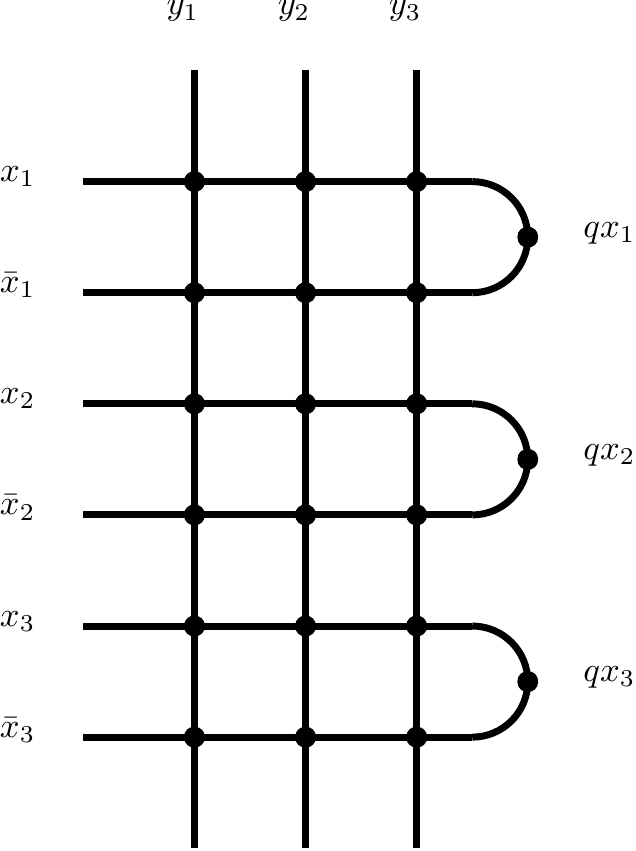}}
\caption{\label{grid} The grid corresponding to VSASMs.}
\end{figure}

\begin{figure}
\scalebox{0.8}{\includegraphics{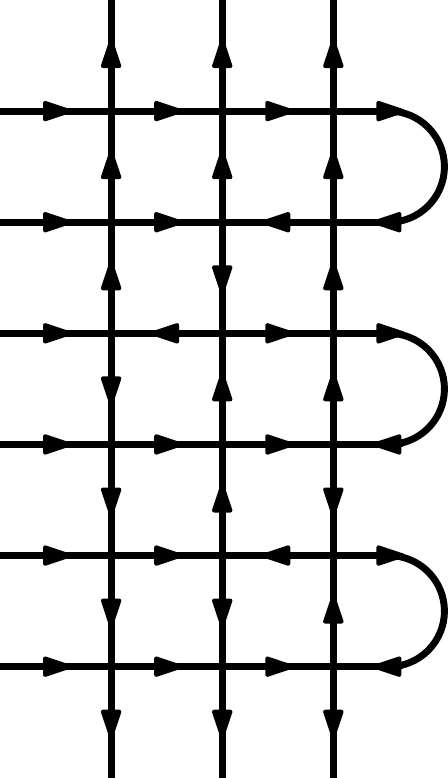}}
\caption{\label{ex} The six-vertex configuration of our example.}
\end{figure}

Next we use the well-known correspondence between ASMs and the six-vertex model as explained in Section~\ref{asm} to translate the $2n \times (n+1)$ matrices into orientations of graphs as indicated in Figure \ref{grid}. In our example we obtain Figure~\ref{ex}. Note that the right boundary, i.e., the fixed column $(1,-1,1,\ldots,-1)^T$, is modeled via U-turns with down-pointing orientation.

\begin{figure}
\begin{center}
\begin{tabular}{lcc}
U-turn & \raisebox{-0.6cm}{\scalebox{0.7}{\includegraphics{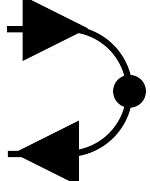}}} & 
\raisebox{-0.6cm}{\scalebox{0.7}{\includegraphics{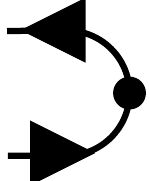}}} \\
weight & $\sigma(b u)$ & $\sigma(b \bar u)$ 
\end{tabular}
\end{center}
\caption{\label{Uweights} Weights of the U-turns.}
\end{figure}

As for the related partition function $Z_U(n; x_1,\ldots,x_n;y_1,\ldots,y_n)$, we allow for the moment both up-pointing and down-pointing U-turns, and the labels and weights are as indicated in Figures~~\ref{grid} and \ref{Uweights}, involving now an additional global parameter $b$. Later we will specialize $b$ so that the weight of a configuration that has at least one up-pointing U-turn other than the top U-turn is $0$. Tsuchiya \cite{Uturn} was the first who derived a formula for this partition function. Here we use Kuperberg's \cite[Theorem 10]{KuperbergRoof} version (up to some normalization factor).

\begin{theo} The U-turn partition function of order $n$ is 
\begin{multline*}
Z_U(n;x_1,\ldots,x_n;y_1,\ldots,y_n) \\
= \frac{\sigma(q^2)^{n- 2 n^2} \prod_{i=1}^{n} \sigma(b \bar y_i) \sigma(q^2 x_i^2) 
\prod_{i,j=1}^{n} \alpha(x_i \bar y_j) \alpha(x_i y_j)}{\prod_{1 \le i < j \le n} \sigma( \bar x_i x_j) \sigma(y_i \bar y_j) \prod_{1 \le i \le j \le n} \sigma(\bar x_i \bar x_j) \sigma(y_i y_j)} \det_{1 \le i,j \le n} \left( \frac{1}{\alpha(x_i \bar y_j)} - \frac{1}{\alpha(x_i y_j)} \right),  
\end{multline*}
where $\alpha(x) = \sigma(q x) \sigma(q \bar x)$.
\end{theo}

In the following, we will specialize 
$$
(x_1,\ldots,x_n)=(x,1,\ldots,1) \quad \text{and} \quad (y_1,\ldots,y_n)=(1,\ldots,1),$$
as well as 
\begin{equation}
\label{special}
b = q \quad \text{and} \quad q + \bar q =1,
\end{equation}
in the partition function. Next we explore how this specialization can be expressed in terms of the numbers $\av(2n+1,i)$. 

First of all, we note that $b = q$ and $x_i = 1$ for $i > 1$ implies that the configurations that have at least one up-pointing U-turn in positions $2,3,\ldots,n$ have weight $0$ and can therefore be omitted. For the remaining configurations we can distinguish between the cases where the topmost U-turn is down-pointing (Case~1) or not (Case~2). 

{\bf Case 1.} If the topmost U-turn is down-pointing, then the top row is forced and all vertex configurations are of type $\ne$. In the second row, there is precisely one configuration of type $\one$, say in position $i$ counted from the left,  and the configurations right of it are all of type $\nw$, while the configurations left of it are of type $\ne$. Such configurations correspond to $(2n+1) \times (2n+1)$ VSASMs that have the first $1$ in the second row in the $i$-th column. The top U-turn contributes $\sigma(q^2x)$, while all other $n-1$ U-turns contribute $\sigma( q^2)$. In total such a configuration has the following weight
$$
\left( \frac{\sigma(q x)}{\sigma(q^2)} \right)^{2n-i} \left( \frac{\sigma(q \bar x)}{\sigma(q^2)} \right)^{i-1} \sigma(q^2x) \sigma( q^2)^{n-1}.
$$
This case contributes the following term towards the partition function
$$
 \sum_{i=1}^{n} \av(2n+1,i) \left( \frac{\sigma(q \bar x)}{\sigma(q x)} \right)^{i} 
\left( \frac{\sigma(q x)}{\sigma(q^2)} \right)^{2 n} \left( \frac{\sigma(q \bar x)}{\sigma(q^2)} \right)^{-1} \sigma(q^2x) \sigma( q^2)^{n-1}.
$$

{\bf Case 2.} If the topmost U-turn is up-pointing, there is a unique occurrence of $\one$ in the top row, say in position $i$. There is either one occurrence of $\one$ in the second row, say in position $j$ with $1 \le j < i$, or no such occurrence. 

In the first case, the weight is 
$$
\left( \frac{\sigma(q x)}{\sigma(q^2)} \right)^{2i-j-2} \left( \frac{\sigma(q \bar x)}{\sigma(q^2)} \right)^{2n-2i+j-1} \sigma(\bar x) \sigma( q^2)^{n-1}.
$$
We notice that for fixed $i$, these configurations give rise to all the configurations counted by $\av(2n+1,j)$. So, this case contributes the following term towards the partition function
$$
\sum_{j=1}^{n} \av(2n+1,j) \left( \frac{\sigma(q \bar x)}{\sigma(q x)} \right)^{j} 
\sum_{i=j+1}^{n} \left( \frac{\sigma(q x)}{\sigma(q^2)} \right)^{2i-2} 
\left( \frac{\sigma(q \bar x)}{\sigma(q^2)} \right)^{2n-2i-2} \sigma( \bar x) \sigma( q^2)^{n-1}.
$$
In the second case the weight is 
$$
\left( \frac{\sigma(q x)}{\sigma(q^2)} \right)^{i} \left( \frac{\sigma(q \bar x)}{\sigma(q^2)} \right)^{2n-i-1} \sigma( \bar x) \sigma( q^2)^{n-1}.
$$ We notice that such configurations are essentially the same as the ones counted by $\av(2n+1,i)$ with just the first U-turn reversed, so this contributes the following term towards the partition function
$$
\sum_{i=1}^{n} \av(2n+1,i) \left( \frac{\sigma(q x)}{\sigma(q \bar x)} \right)^{i} 
\left( \frac{\sigma(q \bar x)}{\sigma(q^2)} \right)^{2n-1} \sigma( \bar x) \sigma( q^2)^{n-1}.
$$

From this it follows that 
\begin{multline}\label{vs-eq-1}
Z_U(n;x,\underbrace{1,\ldots,1}_{n-1};\underbrace{1,\ldots,1}_{n}) 
=  \sum_{i=1}^{n} \av(2n+1,i) \left( \frac{\sigma(q x)}{\sigma(q \bar x)} \right)^{i} 
\left( \frac{\sigma(q \bar x)}{\sigma(q^2)} \right)^{2n-1} \sigma( \bar x) \sigma( q^2)^{n-1}\\
+ \sum_{j=1}^{n} \av(2n+1,j) \left( \frac{\sigma(q \bar x)}{\sigma(q x)} \right)^{j} 
\sum_{i=j+1}^{n} \left( \frac{\sigma(q x)}{\sigma(q^2)} \right)^{2i-2} 
\left( \frac{\sigma(q \bar x)}{\sigma(q^2)} \right)^{2n-2i-2} \sigma( \bar x) \sigma( q^2)^{n-1}\\
+\sum_{i=1}^{n} \av(2n+1,i) \left( \frac{\sigma(q \bar x)}{\sigma(q x)} \right)^{i} 
\left( \frac{\sigma(q x)}{\sigma(q^2)} \right)^{2 n} \left( \frac{\sigma(q \bar x)}{\sigma(q^2)} \right)^{-1} \sigma(q^2x) \sigma( q^2)^{n-1},
\end{multline}
provided equation \eqref{special} is fulfilled. 

In the end, we want to perform the following transformation of variable
$$
z = \frac{\sigma(q \bar x)}{\sigma(q x)}
$$
and eliminate $x$. A tedious but straight forward computation shows that 
\begin{equation}
\label{straightforward}
-\sigma(q^2)^n \sigma(q \bar x)^{-2n} \frac{1+z}{1-2z} Z_U(n;x,\underbrace{1,\ldots,1}_{n-1};\underbrace{1,\ldots,1}_{n}) = 
\sum_{i=1}^{n} \av(2n+1,i) \left(z^{i-2n-1}+z^{-i} \right). 
\end{equation}

We define
$$
W^{-}(\alpha_1,\ldots,\alpha_n;x_1,\ldots,x_n) = \det_{1 \le i,j \le n} 
\left( x_i^{\alpha_j} - x_i^{-\alpha_j} \right)
$$
and 
$$
Sp_{2n}(\lambda_1,\ldots,\lambda_n;x_1,\ldots,x_n) = \frac{W^{-}(\lambda_1+n,\lambda_2+n-1,\ldots,\lambda_n+1;x_1,\ldots,x_n)}{W^{-}(n,n-1,\ldots,1;x_1,\ldots,x_n)}.
$$
Then $Sp_{2n}(\lambda_1,\ldots,\lambda_n;x_1,\ldots,x_n)$ is the character of the irreducible representation of the symplectic group $Sp_{2n}(\mathbb{C})$ corresponding to the partition $(\lambda_1,\ldots,\lambda_n)$. 
Okada \cite[Theorem 2.4]{OkadaCharacters} showed that 
\begin{multline*}
\sigma(q)^{-2n^2+2n} \sigma(q^2)^{2n^2-n} \prod_{i=1}^{n} \sigma(b \bar y_i)^{-1} \sigma(q^2 x_i^2)^{-1} Z_U(n;x_1,\ldots,x_n;y_1,\ldots,y_n) \\ 
= 3^{-n(n-1)} Sp_{4n}(n-1,n-1,n-2,n-2,\ldots,0,0;x_1^2,\ldots,x_n^2,y_1^2,\ldots,y_n^2), 
\end{multline*}
provided $q+\bar q = 1$. From this, we get the following
$$
Z_U(n;x,1,\ldots,1;1,\ldots,1) = \sigma(q^2 x^2) \sigma(q^2)^{n-1} 3^{-n(n-1)} Sp_{4n}(n-1,n-1,\ldots,0,0;x^2,1,\ldots,1).
$$
Combining this with equation \eqref{straightforward}, we obtain 
\begin{multline}
\label{relation}
-\sigma(q^2)^{2n-1} \sigma(q \bar x)^{-2n} \frac{1+z}{1-2 z} \sigma(q^2 x^2) 
3^{-n(n-1)} Sp_{4n}(n-1,n-1,\ldots,0,0;x^2,1,\ldots,1) \\ 
= 
\sum_{i=1}^{n} \av(2n+1,i) \left(z^{i-2n-1}+z^{-i} \right).
\end{multline}

On the other hand, Okada \cite[Theorem 2.5]{OkadaCharacters} also showed that the partition function for off-diagonally symmetric ASMs\footnote{ASMs which are (diagonally) symmetric with a null diagonal are called off-diagonally symmetric ASMs.} (OSASMs) satisfy
\begin{multline*}
\sigma(q)^{-2n^2+2n} \sigma(q^2)^{2n^2-2n} Z_O(n;x_1,\ldots,x_{2n}) \\ 
= 3^{-n(n-1)} Sp_{4n}(n-1,n-1,n-2,n-2,\ldots,0,0;x_1^2,\ldots,x_{2n}^2), 
\end{multline*}
provided $q+\bar q=1$, and where $Z_O(n;x_1,\ldots,x_{2n})$ is the partition function for OSASMs, which was explicitly calculated by Kuperberg \cite[Theorem 10]{KuperbergRoof} (up to some normalization factor) as
\begin{equation}\label{partition-zo}
    Z_O(n;x_1, x_2, \ldots, x_{2n})=\frac{\sigma(q^2)^{2n-2n^2}\prod_{i,j=1}^{2n}\alpha(x_ix_j)}{\prod_{i,j=1}^{2n}\sigma(\bar x_i x_j)}\pff_{1\leq i, j \leq 2n}\left(\frac{\sigma(\bar x_i x_j)}{\alpha(x_i x_j)}\right)
\end{equation}
for general $q$ and used by Okada to obtain the specialization. Here the {\it Pfaffian} of a triangular array $(a_{i,j})_{1\leq i<j\leq 2n}$ is defined as
\[
\pf_{1\leq i, j\leq 2n}(a_{i,j})=\sum_{\mycom{\pi=\{(i_1,j_1),\ldots, (i_n,j_n)\}}{\pi~\text{a perfect matching of} K_{2n}, i_k<j_k}}\sgn \pi \prod_{k=1}^na_{i_k,j_k},
\]
where $\sgn \pi$ is the sign of the permutation $i_1j_1\ldots i_nj_n$.

Thus, it follows that
\begin{equation}\label{eq-zo}
 3^{-n(n-1)} Sp_{4n}(n-1,n-1,n-2,n-2,\ldots,0,0;x^2,1, \ldots,1)= Z_O(n;x,1,\ldots,1),
\end{equation}
provided $q+\bar q=1$ holds. From Razumov and Stroganov \cite[Equation (24)]{RazumovStroganov}, we know that
\begin{equation}\label{rs-zo}
Z_O(n;x,1,\ldots,1)=\sum_{i=2}^{2n}\ao(2n,i)z^{-i+2}\sigma(q\bar x)^{2n-2}\sigma(q^2)^{-2n+2},
\end{equation}
where $\ao(2n,i)$ is the number of $2n \times 2n$ off-diagonally symmetric ASMs with the unique $1$ in the first row in the $i$-th position, given by
\begin{equation}\label{value}
    \ao(2n,i)=
    \begin{cases}
    0, &\text{if~} i=0,1;\\
    \dfrac{1}{2^{n-1}}\displaystyle\prod_{k=1}^{n-1}\dfrac{(6k-2)!(2k-1)!}{(4k-1)!(4k-2)!}\displaystyle\sum_{k=1}^{i-1}(-1)^{i+k-1}\dfrac{(2n+k-2)!(4n-k-1)!}{(4n-2)!(k-1)!(2n-k)!}, &\text{for}~i\geq 2.
    \end{cases}
\end{equation}
Using equations \eqref{eq-zo} and \eqref{rs-zo}, we obtain
\begin{equation}\label{eq:sp-ao}
    Sp_{4n}(n-1,n-1,n-2,n-2,\ldots,0,0;x^2,1, \ldots,1)=3^{n(n-1)}\left(\frac{\sigma(q\bar x)}{\sigma(q^2)}\right)^{2n-2}\sum_{i=2}^{2n}\ao(2n,i)z^{-i+2}.
\end{equation}
Combining this with \eqref{relation}, we obtain
\begin{equation*}
-\sigma(q^2) \sigma(q \bar x)^{-2} \frac{1+z}{1-2 z} \sigma(q^2 x^2) \sum_{i=2}^{2n}\ao(2n,i)z^{-i+2}= 
\sum_{i=1}^{n} \av(2n+1,i) \left(z^{i-2n-1}+z^{-i} \right).
\end{equation*}
Assuming $q+\bar q=1$, we have $-\dfrac{\sigma(q^2)\sigma(q^2x^2)}{\sigma(q\bar x)^2}=\dfrac{1-2z}{z^2}$, and
\begin{equation}\label{final}
\sum_{i=2}^{2n}\ao(2n,i)\left(z^{-i}+z^{-i+1}\right)= 
\sum_{i=1}^{n} \av(2n+1,i) \left(z^{i-2n-1}+z^{-i} \right)
\end{equation}
follows.
Now, we compare the coefficients of $z^{-i}$ for $1\leq i\leq n$ from equation \eqref{final} to conclude
\begin{equation}\label{cni}
    \av(2n+1,i)=\ao(2n,i)+\ao(2n,i+1).
\end{equation}
Combining equations \eqref{cni} and \eqref{value}, we get the following theorem.

\begin{theo}\label{ilse}
The number of $(2n+1)\times (2n+1)$ VSASMs with the first $1$ in its second row at position $i$ is given by
\begin{equation*}
\frac{(2n+i-2)!(4n-i-1)!}{2^{n-1}(4n-2)!(i-1)!(2n-i)!} \prod_{k=1}^{n-1}\frac{(6k-2)!(2k-1)!}{(4k-1)!(4k-2)!}.
\end{equation*}

\end{theo}

\section{Vertically and Horizontally Symmetric ASMs}\label{sec:vhsasm}

ASMs that are invariant under the reflection in the vertical symmetry axis as well as the horizontal symmetry axis also exist only for odd order. In this section we focus on such matrices. Due to a slight difference between how the matrices are enumerated for order $4n+3$ and order $4n+1$, we deal with the refined enumeration of both the cases in separate subsections below. We assume $n\geq 1$, unless otherwise mentioned. The only VHSASM of order $1$ is the single entry matrix $(1)$ and for order $3$, the following matrix
\[
\begin{pmatrix}
 0&1&0\\
 1&-1&1\\
 0&1&0
\end{pmatrix}.
\]

\subsection{VHSASMs of order \texorpdfstring{$4n+3$}{Lg-1}}\label{sec:vh-3}
First, we consider the case for VHSASMs of order $4n+3$. An example of such a VHSASM of order $15$ is 
$$
\begin{pmatrix}
0&0&0&0&0&0&0&1&0&0&0&0&0&0&0\\
0&0&0&1&0&0&0&-1&0&0&0&1&0&0&0\\
0&0&1&-1&0&0&0&1&0&0&0&-1&1&0&0\\
0&0&0&1&0&0&0&-1&0&0&0&1&0&0&0\\
0&0&0&0&0&0&0&1&0&0&0&0&0&0&0\\
0&1&-1&0&0&1&0&-1&0&1&0&0&-1&1&0\\
0&0&0&0&0&0&0&1&0&0&0&0&0&0&0\\
1&-1&1&-1&1&-1&1&-1&1&-1&1&-1&1&-1&1\\
0&0&0&0&0&0&0&1&0&0&0&0&0&0&0\\
0&1&-1&0&0&1&0&-1&0&1&0&0&-1&1&0\\
0&0&0&0&0&0&0&1&0&0&0&0&0&0&0\\
0&0&0&1&0&0&0&-1&0&0&0&1&0&0&0\\
0&0&1&-1&0&0&0&1&0&0&0&-1&1&0&0\\
0&0&0&1&0&0&0&-1&0&0&0&1&0&0&0\\
0&0&0&0&0&0&0&1&0&0&0&0&0&0&0\\
\end{pmatrix}.
$$
From the discussion in the previous section, it is clear that VHSASMs also have two $1$'s in its second row, and the second row is determined by the position of the first $1$ in this row. The middle row of a VHSASMs is $(1,-1,1,\ldots, -1, 1)$ by the horizontal symmetry.  The aim of this subsection is to give a generating function result for the refined enumeration of order $4n+3$ VHSASMs with respect to the position of the first $1$ in the second row.

 From the preceding paragraph, it is clear that $(4n+3) \times (4n+3)$ VHSASMs correspond to $(2n+1) \times (2n+1)$ matrices with entries in $\{\pm 1,0\}$ that have the following properties.
\begin{enumerate}
\item The non-zero entries alternate in each row and column.
\item The topmost non-zero entry of each column is $1$, except the last column which is always $(-1,1,\ldots, 1, -1)^T$.
\item The first non-zero entry of each row is $1$, except for the last row which is always $(-1,1,\ldots, 1, -1)$.
\end{enumerate}
The $7 \times 7$ matrix with these properties that corresponds to the above VHSASM is
$$
\begin{pmatrix}
0&0&1&0&0&0&-1\\
0&1&-1&0&0&0&1\\
0&0&1&0&0&0&-1\\
0&0&0&0&0&0&1\\
1&-1&0&0&1&0&-1\\
0&0&0&0&0&0&1\\
-1&1&-1&1&-1&1&-1\\
\end{pmatrix}.
$$
(We deleted the last $(2n+1)$ columns and rows as well as the first row and first column.) 
Such a matrix has a unique $1$ in its first row and let $\avh(4n+3,i+1)$ be the number of such matrices that have this unique $1$ in column $i$ (the index $i$ runs from $1$ to $2n$, because the last position has a fixed $-1$). Note that $\avh(4n+3,i)$ is now equal to the number of VHSASMs of order $4n+3$ where the first $1$ in the second row is in position $i$. Due to the horizontal symmetry $\avh(2n+1,1)=0$ for all $n$. Now, if we delete the top two rows from the $(2n+1) \times (2n+1)$ matrix we obtain a $(2n-1) \times (2n+1)$ matrix with properties (1) and (3), but property (2) replaced by the following,
\begin{enumerate}
\item[(2')] The topmost non-zero entry of each column is $1$, except for the last column which is always $(-1,1,\ldots, 1, -1)^T$, and one other column whose topmost non-zero entry is $-1$ (if a non-zero entry exists at all in this column).
\end{enumerate}
In our example, we obtain 
$$
\begin{pmatrix}
0&0&1&0&0&0&-1\\
0&0&0&0&0&0&1\\
1&-1&0&0&1&0&-1\\
0&0&0&0&0&0&1\\
-1&1&-1&1&-1&1&-1\\
\end{pmatrix}.
$$

We let $\evh(4n+3,j)$ denote the number of such $(2n-1) \times (2n+1)$ matrices with the special column in (2') being column $j$. When passing from the 
$(2n+1) \times (2n+1)$ matrix to the $(2n-1) \times (2n+1)$ matrix, the column $j$ is the position of the first $1$ in the second row of the $(2n+1) \times (2n+1)$ matrix if the second row contains two $1$'s and otherwise it is the position of the unique $1$ in the top row. We can deduce the following simple relation between $\avh(4n+3,i)$ and $\evh(4n+3,j)$:
\begin{multline}
    \label{dni}
\avh(4n+3,i+1) = \sum_{j=1}^{i} \evh(4n+3,j) \\
\Leftrightarrow \evh(4n+3,i) = \avh(4n+3,i+1) - \avh(4n+3,i).
\end{multline}
Hence, in order to compute $\avh(4n+3,i)$, it suffices to compute $\evh(4n+3,i)$.

\begin{figure}
\begin{center}
	\scalebox{0.7}{
\includegraphics{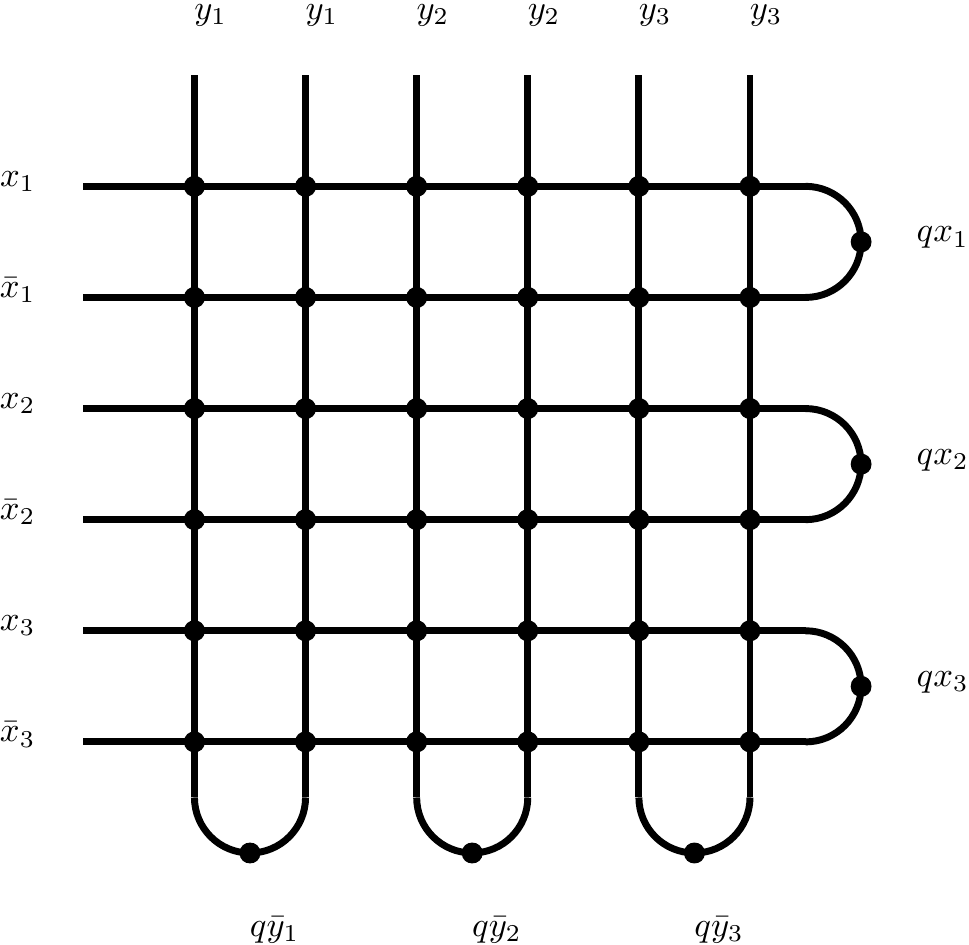}}
\caption{\label{grid-2} The grid corresponding to VHSASMs.}
\end{center}
\end{figure}

\begin{figure}
\scalebox{0.7}{\includegraphics{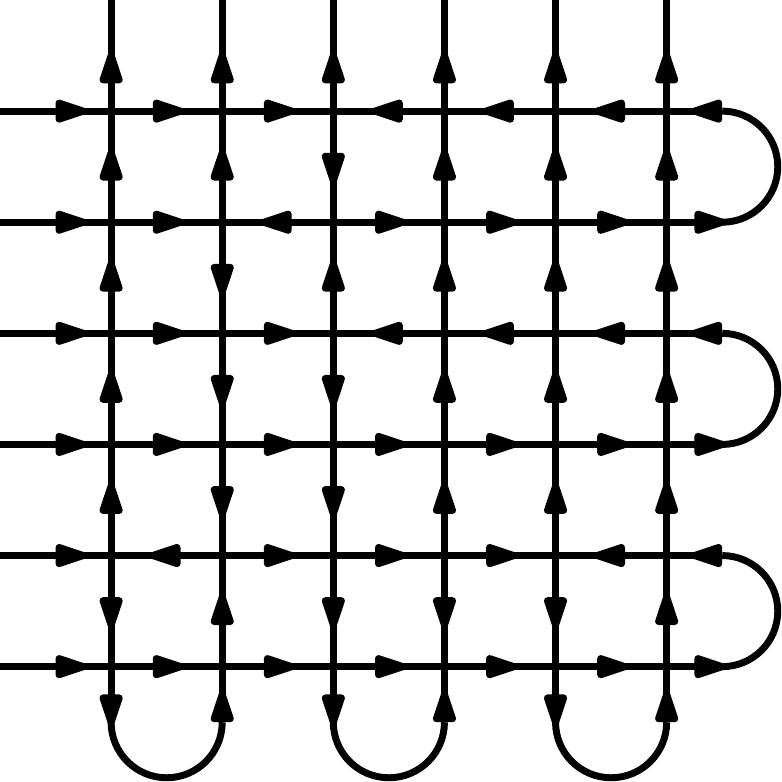}}
\caption{\label{ex-v} The six-vertex configuration of our example.}
\end{figure}

\begin{figure}
\begin{center}
\begin{tabular}{lcccc}
U-turn & \raisebox{-0.6cm}{\scalebox{0.7}{\includegraphics{Udown.pdf}}} & 
\raisebox{-0.6cm}{\scalebox{0.7}{\includegraphics{Uup.pdf}}} & \raisebox{-0.6cm}{\scalebox{0.75}{\includegraphics{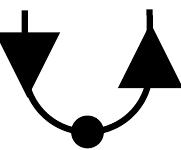}}} & 
\raisebox{-0.6cm}{\scalebox{0.75}{\includegraphics{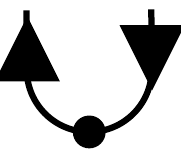}}} \\
weight & $\sigma(b u)$ & $\sigma(b \bar u)$ & $\sigma(c\bar u)$ & $\sigma(c u)$ 
\end{tabular}
\end{center}
\caption{\label{UUweights-2} Weights of the U-turns.}
\end{figure}

Next we use the correspondence between ASMs and the six-vertex model as explained in Section~\ref{asm} to translate the $(2n+1) \times (2n+1)$ matrices into directed graphs (see Figure \ref{grid-2}). In our example we obtain Figure~\ref{ex-v}. Note that the right boundary, i.e., the fixed column $(-1,1,\ldots,-1)^T$, is modeled via U-turns with up-pointing orientation and the bottom boundary, i.e., the fixed row $(-1,1,-1,\ldots, -1)$ is modeled via right pointing U-turns. For the partition function $$Z_{UU}(n;x_1,\ldots,x_n;y_1,\ldots,y_n),$$ we allow both up-pointing and down-pointing U-turns for the right boundary as well as both right-pointing and left-pointing U-turns for the bottom boundary, and the weights are as indicated in Figure~\ref{UUweights-2}, involving now another global parameter $c$.

We shall use the following formula for this partition function that was derived by Kuperberg~\cite[Theorem 10]{KuperbergRoof} (up to some normalization factor).

\begin{theo}\label{thm:vhs-z}
The UU-turn partition function of order $n$ is 
\begin{multline*}
Z_{UU}(n;x_1,\ldots,x_n;y_1,\ldots,y_n) \\
= \frac{\sigma(q^2)^{n- 4 n^2} \prod_{i=1}^{n} \sigma(q^2 \bar y_i^2) \sigma(q^2 x_i^2) 
\prod_{i,j=1}^{n} \alpha(x_i y_j)^2 \alpha(x_i \bar y_j)^2}{\prod_{1 \le i < j \le n} \sigma( \bar x_i x_j)^2 \sigma(y_i \bar y_j)^2 \prod_{1 \le i \le j \le n} \sigma(\bar x_i \bar x_j)^2 \sigma(y_i y_j)^2} \det_{1 \le i, j \le n} M_U \det_{1 \le i, j \le n} M_{UU},  
\end{multline*}
where $\alpha(x) = \sigma(q x) \sigma(q \bar x)$, $$(M_U)_{i,j}= \left( \frac{1}{\alpha(x_i \bar y_j)} - \frac{1}{\alpha(x_i y_j)} \right),$$ and $$(M_{UU})_{i,j}=\ \left( \frac{\sigma(b\bar y_j)\sigma(c x_i)}{\sigma(qx_i \bar y_j)}-\frac{\sigma(b\bar y_j)\sigma(c \bar x_i)}{\sigma(q\bar x_i \bar y_j)} -\frac{\sigma(b y_j)\sigma(c x_i)}{\sigma(qx_i y_j)} +\frac{\sigma(b y_j)\sigma(c \bar x_i)}{\sigma(q\bar x_i  y_j)}  \right),$$
and all determinants are of order n.
\end{theo}

In the following, we will specialize 
$$
(x_1,\ldots,x_n)=(x,1,\ldots,1) \quad \text{and} \quad (y_1,\ldots,y_n)=(1,\ldots,1),$$
as well as 
\begin{equation}
\label{special-2}
b = \bar q, \quad c= \bar q  \quad \text{and} \quad q + \bar q =1,
\end{equation}
in the partition function. First of all, we note that $b = \bar q$ and $x_i = 1$ for $i > 1$ implies that the configurations that have at least one down-pointing U-turn in positions $2,3,\ldots,n$ have weight $0$ and can therefore be omitted. For the remaining configurations we can distinguish between the cases where the topmost U-turn is down-pointing (Case~1) or not (Case~2). Also $c=\bar q$ means that configurations which have at least one left pointing U-turns in the bottom boundary have weight $0$, and so they are omitted.

{\bf Case 1.} If the topmost U-turn is down-pointing, then the top row is forced and all vertex configurations are of type $\ne$. In the second row, there is precisely one configuration of type $\one$, say in position $i$ counted from the left,  and the configurations right of it are all of type $\nw$, while the configurations left of it are of type $\ne$. The top U-turn contributes $\sigma(x)$, while all other $2n-1$ U-turns contribute $\sigma(\bar q^2)$. In total such a configuration has the following weight
$$
\left( \frac{\sigma(q x)}{\sigma(q^2)} \right)^{4n-i} \left( \frac{\sigma(q \bar x)}{\sigma(q^2)} \right)^{i-1} \sigma(x) \sigma(\bar q^2)^{2n-1}.
$$

{\bf Case 2.} In this case, there is a unique occurrence of $\one$ in the top row, say in position $i$. There is either one occurrence of $\one$ in the second row, say in position $j$ with $1 \le j < i$, or no such occurrence. In the first case, the weight is 
$$
\left( \frac{\sigma(q x)}{\sigma(q^2)} \right)^{2i-j-2} \left( \frac{\sigma(q \bar x)}{\sigma(q^2)} \right)^{4n-2i+j-1} \sigma(\bar q^2 \bar x) \sigma(\bar q^2)^{2n-1},
$$
while in the second case the weight is 
$$
\left( \frac{\sigma(q x)}{\sigma(q^2)} \right)^{i} \left( \frac{\sigma(q \bar x)}{\sigma(q^2)} \right)^{4n-i-1} \sigma(\bar q^2 \bar x) \sigma(\bar q^2)^{2n-1}.
$$
From this, it follows that
\begin{multline}\label{vsh-eq-1}
Z_{UU}(n;x,\underbrace{1,\ldots,1}_{n-1};\underbrace{1,\ldots,1}_{n}) 
= \sum_{i=1}^{2n}\evh(4n+3,i) \left( \frac{\sigma(q x)}{\sigma(q \bar x)} \right)^{i} 
\left( \frac{\sigma(q \bar x)}{\sigma(q^2)} \right)^{4n-1} \sigma(\bar q^2 \bar x) \sigma(\bar q^2)^{2n-1} \\
+ \sum_{j=1}^{2n} \evh(4n+3,j) \left( \frac{\sigma(q \bar x)}{\sigma(q x)} \right)^{j} 
\sum_{i=j+1}^{2n} \left( \frac{\sigma(q x)}{\sigma(q^2)} \right)^{2i-2} 
\left( \frac{\sigma(q \bar x)}{\sigma(q^2)} \right)^{4n-2i-2} \sigma(\bar q^2 \bar x) \sigma(\bar q^2)^{2n-1} \\
+ \sum_{i=1}^{2n} \evh(4n+3,i) \left( \frac{\sigma(q \bar x)}{\sigma(q x)} \right)^{i} 
\left( \frac{\sigma(q x)}{\sigma(q^2)} \right)^{4 n} \left( \frac{\sigma(q \bar x)}{\sigma(q^2)} \right)^{-1} \sigma(x) \sigma(\bar q^2)^{2n-1}.
\end{multline}
From here, substituting $z=\dfrac{\sigma(q\bar x)}{\sigma(qx)}$, we shall arrive at
\begin{equation}
\label{straightforward-2}
- \sigma(q^2)^{2n} \sigma(q \bar x)^{-4n} \frac{1-z^2}{1-2z} Z_{UU}(n;x,\underbrace{1,\ldots,1}_{n-1};\underbrace{1,\ldots,1}_{n}) = 
\sum_{i=1}^{2n} \evh(4n+3,i) \left( z^{i-4n-1} - z^{-i+1} \right).
\end{equation}

Okada \cite[Theorem 2.4]{OkadaCharacters} showed that
\begin{multline*}
 \prod_{i=1}^{n} \sigma(q^2 \bar y_i^2)^{-1} \sigma(q^2 x_i^2)^{-1} Z_{UU}(n;x_1,\ldots,x_n;y_1,\ldots,y_n) \\ 
= 3^{-2n^2+n} Sp_{4n}(n-1,n-1,n-2,n-2,\ldots,0,0;x_1^2,\ldots,x_n^2,y_1^2,\ldots,y_n^2)\\ \times Sp_{4n+2}(n,n-1,n-1,n-2,n-2,\ldots,0,0;x_1^2,\ldots,x_n^2,y_1^2,\ldots,y_n^2,1), 
\end{multline*}
provided \eqref{special-2} is satisfied. From this, we get the following
\begin{multline}\label{vsh-2}
Z_{UU}(n;x, 1,\ldots,1;1,\ldots,1)  
=\sigma(q^2)^{2n-1}\sigma(q^2x^2) 3^{-2n^2+n} \\ \times Sp_{4n}(n-1,n-1,n-2,n-2,\ldots,0,0;x^2,1,\ldots,1)\\ \times Sp_{4n+2}(n,n-1,n-1,n-2,n-2,\ldots,0,0;x^2,1,\ldots,1).
\end{multline}
Combining equations \eqref{straightforward-2} and \eqref{vsh-2}, we get
\begin{multline}\label{vsh-3}
-3^{-2n^2+n}\sigma(q^2x^2) \sigma(q^2)^{4n-1} \sigma(q \bar x)^{-4n} \frac{1-z^2}{1-2z} \\ \times Sp_{4n}(n-1,n-1,n-2,n-2,\ldots,0,0;x^2,1,\ldots,1)\\ \times Sp_{4n+2}(n,n-1,n-1,n-2,n-2,\ldots,0,0;x^2,1,\ldots,1)  \\ = 
\sum_{i=1}^{2n} \evh(4n+3,i) \left( z^{i-4n-1} - z^{-i+1} \right). 
\end{multline}

From equation \eqref{gf} in Appendix \ref{appen}, we have 
\begin{equation}\label{sp-4n+2}
 Sp_{4n+2}(n,n-1,n-1,n-2,n-2,\ldots,0,0;x^2,1,\ldots,1) = \sum_{1 \le j \le i \le n+1} Q_{n,i} x^{2i-4j+2},
\end{equation}
where
\begin{multline}\label{value-qni}
Q_{n,i} = \frac{3^{n(n-1)}}{2^{n-1} (4n-1)!} \prod_{j=0}^{n-1}  \frac{(4j+3)(6j+6)!}{(2n+2j+1)!} 
\sum_{j=0}^{n}  \left[ \frac{ 9^j (3j-2n-i+2)_{4n-3} (3n-3j+1)}{(3j)! (n-j)! (3j+1)_{3n}} \right. 
\\
\left.  \times \left( \frac{ \left( n-j + \frac{4}{3} \right)_{2j} (2n+3j-i-1)_2}{(3n+3j+1)_2} -
 \frac{\left( n-j + \frac{2}{3} \right)_{2j} (-2n+3j-i)_2}{(3n-3j+1)_2} \right) \right]
\end{multline}
with $(a)_n=a(a+1) \cdots (a+n-1)$. The purpose of Appendix \ref{appen} is also to provide a combinatorial interpretation of $Q_{n,i}$ in terms of rhombus tilings.

Using equations \eqref{eq:sp-ao}, \eqref{vsh-3}, \eqref{sp-4n+2} 
as well as $q+\bar q=1$ and $$\evh(4n+3,i)=\avh(4n+3,i+1)-\avh(4n+3,i)$$ with some simplifications, we now arrive at the following equation
\begin{multline}\label{vhsasm-apr-1}
3^{-n^2} (z^2 - z+1)^n (1-z^2) \left(\sum_{i=2}^{2n}\ao(2n,i)z^{-i}\right)\left( \sum_{1 \le j \le i \le n+1} Q_{n,i} x^{2i-4j+2} \right) \\= 
\sum_{i=1}^{2n+1}  \left( \avh(4n+3,i+1) - \avh(4n+3,i) \right) \left( z^{i-2n-1} - z^{-i+2n+1} \right)  .
\end{multline}
Replacing $x^2=\dfrac{zq-1}{q-z}$, using $q+\bar q=1$ we get 
\begin{multline*}
3^{-n^2} (1-z^2) \left(\sum_{i=2}^{2n}\ao(2n,i)z^{-i}\right) \left( \sum_{1 \le j \le i \le n+1} Q_{n,i} (z q -1)^{n+i-2j+1} (q-z)^{n-i+2j-1}  (-q)^{-n}  \right)\\= 
\sum_{i=1}^{2n+1}  \left( \avh(4n+3,i+1) - \avh(4n+3,i) \right) \left( z^{i-2n-1} - z^{-i+2n+1} \right)  .
\end{multline*}

Thus, we have proved the following result.

\begin{theo}\label{thm:vh}
Let $\avh(4n+3,i)$ denote the number of VHSASMs of order $4n+3$, with the first occurrence of a $1$ in the second row be in the $i$-th column. Then,  for all $n\geq 1$ the following is satisfied
\begin{multline*}
3^{-n^2}(1-z^2) \left(\sum_{i=2}^{2n}\ao(2n,i)z^{-i}\right) \left( \sum_{1 \le j \le i \le n+1} Q_{n,i} (z q -1)^{n+i-2j+1} (q-z)^{n-i+2j-1} (-q)^{-n} \right)\\= 
\sum_{i=1}^{2n+1}  \left( \avh(4n+3,i+1) - \avh(4n+3,i) \right) \left( z^{i-2n-1} - z^{-i+2n+1} \right),
\end{multline*}
where every quantity appearing on the left-hand side is explicitly known, and $\avh(4n+3,1)=0$ for all $n$.
\end{theo}

\begin{rem}
We can write equation \eqref{vsh-3} differently by using results of Ayyer and Behrend \cite{AyyerBehrend}\footnote{First use equation (7), then Proposition 5, followed by equation (8) and finally Corollary 11, equation (55).} as
\begin{multline*}
-3^{-2n^2-3n-2}(z-1-z^2)^{2n+1}(1-z^2) s(2n+1,2n,2n,\ldots,1,1,0,0;x^2,\bar x^2,1, \ldots,1)\\
=\sum_{i=1}^{2n} \evh(4n+3,i) \left( z^{i-1} - z^{4n-i+1} \right). 
\end{multline*}
where
$$
s(\lambda_1,\ldots,\lambda_n;x_1,\ldots,x_n)=\frac{V(\lambda_1+n-1,\lambda_2+n-2,\ldots,\lambda_n;x_1,\ldots,x_n)}{V(n-1,n-2\ldots,1,0;x_1,\ldots,x_n)}
$$
is the character of the irreducible representation of the general linear group $GL_n(\mathbb{C})$ corresponding to the partition $(\lambda_1,\ldots,\lambda_n)$ (i.e., a Schur function) and $$V(\alpha_1, \ldots, \alpha_n;x_1, \ldots, x_n)=\det_{1\leq i,j\leq n}\left(x_i^{\alpha_j}\right).$$ From this and using the relation \eqref{dni}, we shall arrive that
\begin{multline}\label{schur-vhs-3}
-3^{-2n^2-3n-1}(z-1-z^2)^{2n+1}(1-z^2) s(2n+1,2n,2n,\ldots,1,1,0,0;x^2,\bar x^2,1, \ldots,1)\\
=\sum_{i=1}^{2n} (\avh(4n+3,i+1)-\avh(4n+3,i)) \left( z^{i+1} - z^{4n-i+3} \right). 
\end{multline}
\end{rem}

\subsection{VHSASMs of order \texorpdfstring{$4n+1$}{Lg-2}}\label{subsec:vhs-1}

We now focus on the VHSASMs of order $4n+1$. An example of such a matrix of order $9$ is
\[
\begin{pmatrix}
\textcolor{red}{0}&\textcolor{red}{0}&\textcolor{red}{0}&\textcolor{red}{0}&\textcolor{red}{1}&\textcolor{red}{0}&\textcolor{red}{0}&\textcolor{red}{0}&\textcolor{red}{0}\\
0&\textcolor{red}{0}&\textcolor{red}{0}&\textcolor{red}{1}&\textcolor{red}{-1}&\textcolor{red}{1}&\textcolor{red}{0}&\textcolor{red}{0}&0\\
0&0&\textcolor{red}{0}&\textcolor{red}{0}&\textcolor{red}{1}&\textcolor{red}{0}&\textcolor{red}{0}&0&0\\
0&1&0&\textcolor{red}{0}&\textcolor{red}{-1}&\textcolor{red}{0}&0&1&0\\
1&-1&1&-1&\textcolor{red}{1}&-1&1&-1&1\\
0&1&0&0&-1&0&0&1&0\\
0&0&0&0&1&0&0&0&0\\
0&0&0&1&-1&1&0&0&0\\
0&0&0&0&1&0&0&0&0
\end{pmatrix}.
\]
(The colour is to be ignored for this section. This will be used in subsequent sections.)

The observations on the two top rows of order $4n+3$ VHSASMs follow in this case as well. The aim of this subsection is to give a generating function result for the refined enumeration of VHSASMs with respect to the first $1$ in the second row. However, we need to modify our arguments for this case, as the six vertex configurations of VHSASMs of order $4n+1$ are slightly different than for order $4n+3$.

It is clear that any order $4n+1$ VHSASM corresponds to a $(2n+1)\times (2n+1)$ matrix with the following properties.
\begin{enumerate}
\item The non-zero entries alternate in each row and column.
\item The topmost non-zero entry of each column is $1$; the last column is equal to $(1,-1,1,\ldots,-1)^T$.
\item The first non-zero entry of each row is $1$; the last row is equal to $(1,-1,1,\ldots,-1)$.
\end{enumerate}
The $5\times 5$ matrix with these properties that corresponds to the VHSASM from above is
$$
\begin{pmatrix}
0&0&0&0&1\\
0&0&0&1&-1\\
0&0&0&0&1\\
0&1&0&0&-1\\
1&-1&1&-1&1
\end{pmatrix}.$$
(We deleted the last $2n$ columns and rows.) Let us denote the number of order $4n+1$ VHSASMs with the first occurrence of a $1$ in its second row being placed in the $i$-th column to be $\avh(4n+1,i)$. That is, $\avh(4n+1,i)$ counts the number of $(2n+1)\times (2n+1)$ matrices described above where the unique $1$ in its second row is situated in the $i$-th column. 

Throughout the remainder of this section, we consider $n>1$. It is easy to see that for $n=1$, $\avh(5,1)=0$ and $\avh(5,2)=1$.

We shall use the correspondence between ASMs and the six-vertex model as in the previous section. The grid for order $4n+1$ VHSASMs is the same as in Figure \ref{grid-2}. However, for order $4n+1$ VHSASMs the U-turns in the right boundary are now down-pointing, as opposed to the up-pointing ones for order $4n+3$ VHSASMs and the U-turns on the bottom are now left-pointing, as opposed to the right-pointing ones for order $4n+3$ VHSASMs  (see Figure \ref{ex-v}). This is because we do not delete the first row and the first column here. Again, for the partition function $Z_{UU}(n;x_1,\ldots,x_n;y_1,\ldots,y_n)$, we allow both up-pointing and down-pointing U-turns for the right boundary as well as both right-pointing and left-pointing U-turns for the bottom boundary, and the weights are as indicated in Figure~\ref{UUweights-2}. The partition function is still the one given in Theorem \ref{thm:vhs-z}.

In the following, we will specialize 
$$
(x_1,\ldots,x_n)=(x,1,\ldots,1) \quad \text{and} \quad (y_1,\ldots,y_n)=(1,\ldots,1),$$
as well as 
\begin{equation}
\label{special-3}
b = q, \quad c=  q  \quad \text{and} \quad q + \bar q =1,
\end{equation}
in the partition function. This will give us two cases similar to the cases described in the previous section on VSASMs. Analogous to equation \eqref{vs-eq-1} we shall get the following
\begin{multline}\label{vhsasm-1}
Z_{UU}(n;x,\underbrace{1,\ldots,1}_{n-1};\underbrace{1,\ldots,1}_{n}) 
= \sum_{i=1}^{2n} \avh(4n+1,i) \left( \frac{\sigma(q x)}{\sigma(q \bar x)} \right)^{i} 
\left( \frac{\sigma(q \bar x)}{\sigma(q^2)} \right)^{4n-1} \sigma( \bar x) \sigma( q^2)^{2n-1}  \\
+ \sum_{j=1}^{2n} \avh(4n+1,j) \left( \frac{\sigma(q \bar x)}{\sigma(q x)} \right)^{j} 
\sum_{i=j+1}^{2n} \left( \frac{\sigma(q x)}{\sigma(q^2)} \right)^{2i-2} 
\left( \frac{\sigma(q \bar x)}{\sigma(q^2)} \right)^{4n-2i-2} \sigma( \bar x) \sigma( q^2)^{2n-1}\\
+ \sum_{i=1}^{2n} \avh(4n+1,i) \left( \frac{\sigma(q \bar x)}{\sigma(q x)} \right)^{i} 
\left( \frac{\sigma(q x)}{\sigma(q^2)} \right)^{4 n} \left( \frac{\sigma(q \bar x)}{\sigma(q^2)} \right)^{-1} \sigma(q^2x) \sigma( q^2)^{2n-1},
\end{multline}
provided \eqref{special-3} is satisfied.

Like earlier, we want to perform the following transformation of variable
$$
z = \frac{\sigma(q \bar x)}{\sigma(q x)}
$$
and eliminate $x$. After a straightforward calculation, analogous to how we obtained equation \eqref{straightforward}, we shall get the following
\begin{equation}
\label{straightforward-vhs}
-\sigma(q^2)^{2n} \sigma(q \bar x)^{-4n} \frac{1+z}{1-2z} Z_{UU}(n;x,\underbrace{1,\ldots,1}_{n-1};\underbrace{1,\ldots,1}_{n}) = 
\sum_{i=1}^{2n} \avh(4n+1,i) \left( z^{i-4n-1} + z^{-i} \right). 
\end{equation}

We define
$$
W^{+}(\alpha_1,\ldots,\alpha_n;x_1,\ldots,x_n) = \det_{1 \le i,j \le n} 
\left( x_i^{\alpha_j} + x_i^{-\alpha_j} \right)
$$
and 
$$
O_{2n}(\lambda_1,\ldots,\lambda_n;x_1,\ldots,x_n) = \frac{2W^{+}(\lambda_1+n-1,\lambda_2+n-2,\ldots,\lambda_n;x_1,\ldots,x_n)}{W^{+}(n-1,n-2\ldots,0;x_1,\ldots,x_n)}.
$$
Then $O_{2n}(\lambda_1,\ldots,\lambda_n;x_1,\ldots,x_n)$ is the character of the irreducible representation of the double cover of the even orthogonal group $O_{2n}$ corresponding to the partition $(\lambda_1,\ldots,\lambda_n)$, if $\lambda_n\neq 0$.\footnote{It suffices for our purposes to have only this case.} Okada \cite[Theorem 2.4]{OkadaCharacters} showed that 
\begin{multline*}
\prod_{i=1}^{n}(x_i+\bar x_i)(y_i+\bar y_i)\prod_{i=1}^{n} \sigma(q^2 \bar y_i^2)^{-1} \sigma(q^2 x_i^2)^{-1} Z_{UU}(n;x_1,\ldots,x_n;y_1,\ldots,y_n) \\ 
= 3^{-2n^2+n} Sp_{4n}(n-1,n-1,n-2,n-2,\ldots,0,0;x_1^2,\ldots,x_n^2,y_1^2,\ldots,y_n^2)\\ \times O_{4n}\left(n+\frac{1}{2},n-\frac{1}{2},n-\frac{1}{2},n-\frac{3}{2},n-\frac{3}{2},\ldots,\frac{3}{2},\frac{3}{2},\frac{1}{2};x_1^2,\ldots,x_n^2,y_1^2,\ldots,y_n^2 \right), 
\end{multline*}
provided \eqref{special-3} is satisfied.

The above for our special values gives us
\begin{multline}\label{even-ortho}
2^{2n-1}(x+\bar x) Z_{UU}(n;x,\underbrace{1,\ldots,1}_{n-1};\underbrace{1,\ldots,1}_{n})  \\ 
= \sigma(q^2)^{2n-1}\sigma(q^2x^2)3^{-2n^2+n} Sp_{4n}(n-1,n-1,n-2,n-2,\ldots,0,0;x^2,1,\ldots,1)\\ \times O_{4n}\left(n+\frac{1}{2},n-\frac{1}{2},n-\frac{1}{2},n-\frac{3}{2},n-\frac{3}{2},\ldots,\frac{3}{2},\frac{3}{2},\frac{1}{2};x^2,1,\ldots,1\right).
\end{multline}
Now, by using a formula of Ayyer and Behrend \cite[Proposition 5, then use equation (7)]{AyyerBehrend} we can rewrite equation \eqref{even-ortho} as follows.
\begin{multline}\label{sp-ortho}
 Z_{UU}(n;x,\underbrace{1,\ldots,1}_{n-1};\underbrace{1,\ldots,1}_{n}) 
= 3^{-2n^2+3n-2}\sigma(q^2)^{2n-1}\sigma(q^2x^2)\\
\times (x^2+1+\bar x^2)Sp_{4n}(n-1,n-1,n-2,n-2,\ldots,0,0;x^2,1,\ldots,1)\\ \times Sp_{4n-2}(n-1,n-2,n-2,\ldots,1,0;x^2,1,\ldots,1)
\end{multline}

From equations \eqref{eq:sp-ao}, \eqref{sp-4n+2}, \eqref{straightforward-vhs} and \eqref{sp-ortho}, after some simplifications we shall get the following equation:
\begin{multline*}
    (-1)^{n+1}3^{-n^2+2n-1}(1+z)(z-1-z^2)^{n-1}\left(\sum_{i=2}^{2n}\ao(2n,i)z^{-i}\right) \left(\sum_{1\leq i\leq j \leq n}Q_{n-1,i}x^{2i-4j+2}\right) \\ =\sum_{i=1}^{2n}\avh(4n+1,i)(z^{i-2n-2}+z^{2n-1-i}),
\end{multline*}
where $Q_{n,i}$ is given by \eqref{value-qni}. From the above, using $x^2=\dfrac{zq-1}{q-z}$ we shall arrive at the following result.
\begin{theo}\label{thm:vh-1}
Let $\avh(4n+1,i)$ denote the number of VHSASMs of order $4n+1$, with the first occurrence of a $1$ in the second row be in the $i$-th column. Then, for all $n>1$ the following is satisfied
\begin{multline*}
   3^{-n^2+2n-1} (1+z)\left(\sum_{i=2}^{2n}\ao(2n,i)z^{-i}\right)
    \left(\sum_{1\leq i\leq j \leq n}Q_{n-1,i}(zq-1)^{n+i-2j}(q-z)^{n-i+2j-2} (-q)^{-n+1} \right)\\=\sum_{i=1}^{2n}\avh(4n+1,i)(z^{i-2n-2}+z^{2n-1-i}),
\end{multline*}
where every quantity appearing in the left-hand side is explicitly known.
\end{theo}

\begin{rem}
We can write equation \eqref{even-ortho} differently by using a result of Ayyer and Behrend \cite[Equation (8) and then use Corollary 11, equation (54)]{AyyerBehrend} as
\begin{multline*}
Z_{UU}(n;x,\underbrace{1,\ldots,1}_{n-1};\underbrace{1,\ldots,1}_{n}) 
= (-1)^{n^2}3^{-2n^2-n}\sigma(q^2)^{2n-1}\sigma(q^2x^2)(z^2+1-z)z^{-1}\\ \times s(2n,2n-1,2n-1,\ldots,1,1,0,0;x^2,\bar x^2,1, \ldots,1).
\end{multline*}
From this and using equation \eqref{straightforward-vhs}, we shall arrive that
\begin{multline}\label{schur-vh-1}
    (-1)^{n^2}3^{-2n^2-n}(z-1-z^2)^{2n}(1+z)
 s(2n,2n-1,2n-1,\ldots,1,1,0,0;x^2,\bar x^2,\ldots, 1)
\\= \sum_{i=1}^{2n} \avh(4n+1,i) \left(z^{i}+z^{4n+1-i} \right). 
\end{multline}
\end{rem}

\section{Vertically and Horizontally Perverse ASMs}\label{sec:vhp}

A $(4n+1)\times (4n+3)$ matrix is called a vertically and horizontally perverse ASM (VHPASM) if it satisfies the alternating sign conditions and has the same symmetries as a VHSASM, except the central entry ($\star$) which has opposite signs when read horizontally and vertically. An example of such a VHPASM of dimension $9\times 11$ is
$$
\begin{pmatrix}
0&0&0&0&0&1&0&0&0&0&0\\
0&0&1&0&0&-1&0&0&1&0&0\\
0&0&0&0&0&1&0&0&0&0&0\\
0&1&-1&1&0&-1&0&1&-1&1&0\\
1&-1&1&-1&1&\star&1&-1&1&-1&1\\
0&1&-1&1&0&-1&0&1&-1&1&0\\
0&0&0&0&0&1&0&0&0&0&0\\
0&0&1&0&0&-1&0&0&1&0&0\\
0&0&0&0&0&1&0&0&0&0&0
\end{pmatrix}.
$$
This class was first considered by Kuperberg \cite{KuperbergRoof}.

A VHPASM with the dimensions stated above is said to be of order $4n+2$. It is clear that we can ask for the refined enumerations of VHPASMs with respect to the position of the first occurrence of a $1$ in the second column as well as in the second row\footnote{These numbers are different due to the different lengths of the rows and columns.}. Let us denote by $\avhpc(4n+2,i)$ (resp. $\avhpr(4n+2,i)$) the number of order $4n+2$ VHPASMs with the first occurrence of a $1$ in the second column (resp. row) at the $i$-th row (resp. $i$-th column). In this section, we give enumeration results for these numbers. Since, the technique is similar to the one used in Sections \ref{sec:vsasm} and \ref{sec:vhsasm}, for the sake of brevity we omit certain easily verifiable details. We also assume $n\geq 1$ unless otherwise mentioned. For $n=0$, there are no VHPASMs because of the restriction imposed by the entry $\star$.

It is clear that VHPASMs of order $4n+2$ correspond to $(2n+1) \times (2n+1)$ matrices with entries in $\{\pm 1,0, \star\}$ that have the following properties.
\begin{enumerate}
\item The non-zero entries alternate in each row and column.
\item The topmost non-zero entry of each column is $1$; the last column is $(1,-1,\ldots,-1,\star)^T$. 
\item The first non-zero entry of each row is $1$; the last row is $(-1,1,\ldots, 1, \star)$.
\end{enumerate}
The $5\times 5$ matrix with these properties that corresponds to the VHPASM from above is
$$
\begin{pmatrix}
0&0&0&0&1\\
0&1&0&0&-1\\
0&0&0&0&1\\
1&-1&1&0&-1\\
-1&1&-1&1&\star\\
\end{pmatrix}.
$$
(We deleted the first column, the last $2n+1$ columns as well as the bottom $2n$ rows.)

We again use the correspondence between ASMs and the six-vertex model. The grid for order $4n+2$ VHPASMs is the same as in Figure \ref{grid-2} (we ignore the special entry $\star$ to get the grid). However, for order $4n+2$ VHPASMs the U-turns in the right boundary are now down-pointing and the U-turns on the bottom boundary are right-pointing. Again, for the partition function $Z_{UU}(n;x_1,\ldots,x_n;y_1,\ldots,y_n)$, we allow both up-pointing and down-pointing U-turns for the right boundary as well as both right-pointing and left-pointing U-turns for the bottom boundary, and the weights are as indicated in Figures \ref{weights} and \ref{UUweights-2}. The partition function is still the one given in Theorem \ref{thm:vhs-z}.

In the following, we will specialize 
\[
(x_1,\ldots,x_n)=(x,1,\ldots,1) \quad \text{and} \quad (y_1,\ldots,y_n)=(1,\ldots,1),
\]
as well as 
\begin{equation}
\label{special-4}
b = q, \quad c=  \bar q  \quad \text{and} \quad q + \bar q =1,
\end{equation}
for refined enumeration with respect to rows and
\begin{equation}
\label{special-5}
b = \bar q, \quad c=  q  \quad \text{and} \quad q + \bar q =1,
\end{equation}
for refined enumeration with respect to columns, in the partition function. This will give us two cases similar to the cases described in the previous sections. Analogous to equations \eqref{vs-eq-1} and \eqref{vsh-eq-1}, we shall get the following sets of equations.
\begin{multline}\label{vhp-1}
(-1)^nZ_{UU}(n;x,\underbrace{1,\ldots,1}_{n-1};\underbrace{1,\ldots,1}_{n}) \\
=\sum_{j=1}^{2n} \avhpr(4n+2,j+1) \left( \frac{\sigma(q \bar x)}{\sigma(q x)} \right)^{j} 
\sum_{i=j+1}^{2n} \left( \frac{\sigma(q x)}{\sigma(q^2)} \right)^{2i-2} 
\left( \frac{\sigma(q \bar x)}{\sigma(q^2)} \right)^{4n-2i-2} \sigma( \bar x) \sigma( q^2)^{2n-1} \\
+\sum_{i=1}^{2n} \avhpr(4n+2,i+1) \left( \frac{\sigma(q \bar x)}{\sigma(q x)} \right)^{i} 
\left( \frac{\sigma(q x)}{\sigma(q^2)} \right)^{4 n} \left( \frac{\sigma(q \bar x)}{\sigma(q^2)} \right)^{-1} \sigma(q^2x) \sigma(q^2)^{2n-1} \\
+ \sum_{i=1}^{2n}\avhpr(4n+2,i+1) \left( \frac{\sigma(q x)}{\sigma(q \bar x)} \right)^{i} 
\left( \frac{\sigma(q \bar x)}{\sigma(q^2)} \right)^{4n-1} \sigma(\bar x) \sigma( q^2)^{2n-1}
\end{multline}
provided \eqref{special-4} holds, and \begin{multline}\label{vhp-2}
(-1)^nZ_{UU}(n;x,\underbrace{1,\ldots,1}_{n-1};\underbrace{1,\ldots,1}_{n}) \\
= \sum_{j=1}^{2n} \evhp(4n+2,j) \left( \frac{\sigma(q \bar x)}{\sigma(q x)} \right)^{j} 
\sum_{i=j+1}^{2n} \left( \frac{\sigma(q x)}{\sigma(q^2)} \right)^{2i-2} 
\left( \frac{\sigma(q \bar x)}{\sigma(q^2)} \right)^{4n-2i-2} \sigma(\bar q^2 \bar x) \sigma(\bar q^2)^{2n-1} \\
+\sum_{i=1}^{2n} \evhp(4n+2,i) \left( \frac{\sigma(q \bar x)}{\sigma(q x)} \right)^{i} 
\left( \frac{\sigma(q x)}{\sigma(q^2)} \right)^{4 n} \left( \frac{\sigma(q \bar x)}{\sigma(q^2)} \right)^{-1} \sigma(x) \sigma(\bar q^2)^{2n-1} \\
+ \sum_{i=1}^{2n}\evhp(4n+2,i) \left( \frac{\sigma(q x)}{\sigma(q \bar x)} \right)^{i} 
\left( \frac{\sigma(q \bar x)}{\sigma(q^2)} \right)^{4n-1} \sigma(\bar q^2 \bar x) \sigma(\bar q^2)^{2n-1},
\end{multline}
provided \eqref{special-5} holds and where
\[
\evhp(n,i)=\avhpc(n,i)-\avhpc(n,i-1).
\]

From here, substituting $z=\dfrac{\sigma(q\bar x)}{\sigma(qx)}$, equation \eqref{vhp-1} gives us
\begin{multline}
    \label{straightforward-vhp-1}
(-1)^n\sigma(q^2)^{2n} \sigma(q \bar x)^{-4n} \frac{1+z}{2z-1} Z_{UU}(n;x,\underbrace{1,\ldots,1}_{n-1};\underbrace{1,\ldots,1}_{n})\\ = 
\sum_{i=1}^{2n} \avhpr(4n+2,i+1) \left( z^{i-4n-1}+ z^{-i} \right). 
\end{multline}
provided \eqref{special-4} holds, and equation \eqref{vhp-2} gives us
\begin{multline}
    \label{straightforward-vhp-2}
(-1)^{n} \sigma(q^2)^{2n} \sigma(q \bar x)^{-4n} \frac{1-z^2}{2z-1} Z_{UU}(n;x,\underbrace{1,\ldots,1}_{n-1};\underbrace{1,\ldots,1}_{n})\\ = 
\sum_{i=1}^{2n} \evhp(4n+2,i) \left( z^{i-4n-1} - z^{-i+1} \right),
\end{multline}
provided \eqref{special-5} holds.

Okada \cite[Theorem 2.4]{OkadaCharacters} showed that, if \eqref{special-4} holds, then
\begin{multline}\label{okada-vhp}
(-1)^n\prod_{i=1}^{n} \sigma(q^2 \bar y_i^2)^{-1} \sigma(q^2 x_i^2)^{-1}(y_i^2+1+y_i^{-2})^{-1} Z_{UU}(n;x_1,\ldots,x_n;y_1,\ldots,y_n) \\ 
= 3^{-2n^2+n} \left(Sp_{4n}(n-1,n-1,n-2,n-2,\ldots,0,0;x_1^2,\ldots,x_n^2,y_1^2,\ldots,y_n^2)\right)^2,
\end{multline}
and if \eqref{special-5} holds then
\begin{multline}\label{okada-vhp-new}
(-1)^n\prod_{i=1}^{n} \sigma(q^2 \bar y_i^2)^{-1} \sigma(q^2 x_i^2)^{-1}(x_i^2+1+x_i^{-2})^{-1} Z_{UU}(n;x_1,\ldots,x_n;y_1,\ldots,y_n) \\ 
= 3^{-2n^2+n} \left(Sp_{4n}(n-1,n-1,n-2,n-2,\ldots,0,0;x_1^2,\ldots,x_n^2,y_1^2,\ldots,y_n^2)\right)^2,
\end{multline}
Using equations \eqref{eq:sp-ao}, \eqref{straightforward-vhp-1}, \eqref{straightforward-vhp-2}, \eqref{okada-vhp} and \eqref{okada-vhp-new}, after some simplification we shall arrive at
\begin{equation}\label{vhp-1-final}
    (z^3+1)\left(\sum_{i=2}^{2n}\ao(2n,i)z^{-i}\right)^2=\sum_{i=1}^{2n} \avhpr(4n+2,i+1) \left( z^{i-4n-1}+ z^{-i} \right),
\end{equation}
and
\begin{equation}\label{vhp-2-final}
   (1-z^2)\left(\sum_{i=2}^{2n}\ao(2n,i)z^{-i}\right)^2  = 
\sum_{i=1}^{2n} \evhp(4n+2,i) \left( z^{i-4n-1} - z^{-i+1} \right).
\end{equation}
Now, by comparing coefficients in \eqref{vhp-1-final} and \eqref{vhp-2-final} we obtain the following results.

\begin{theo}
The number of order $4n+2$ VHPASMs with the leftmost occurrence of $1$ in the second row in $i$-th column is
\[
\sum_{k=0}^{i-2}\ao(2n,k+2)\left(\ao(2n,i-k)+\ao(2n,i-3-k)\right),
\]
where $\ao(2n,j)$ is given by \eqref{value} and we take $\ao(2n,-1)=0$.
 \end{theo}
\begin{theo}
The number of order $4n+2$ VHPASMs with the topmost occurrence of $1$ in the second column in the $i$-th row is
\[
\sum_{k=0}^{i-2}\ao(2n,k+2)\left(\ao(2n,i-k)+\ao(2n,i-1-k)\right),
\]
where $\ao(2n,j)$ is given by \eqref{value}.
\end{theo}

\begin{rem}
From above it follows that
\begin{equation}\label{vhp-rel}
    \avhpc(4n+2,i)=\avhpr(4n+2,i)+\avhpc(4n+2,i-1)-\avhpc(4n+2,i-2).
\end{equation}
\end{rem}

\section{Off-Diagonally and Off-Antidiagonally Symmetric ASMs}\label{sec:oos}

An ASM which is diagonally and antidiagonally symmetric with each entry in the diagonal and antidiagonal equal to $0$ are called off-diagonally and off-antidiagonally symmetric ASMs (OOSASMs). These matrices occur for order $4n$ and no product formula for their enumeration is currently known or conjectured. However, Ayyer, Behrend and the first author \cite{AyyerBehrendFischer} introduced the concept of an odd order OOSASM while studying extreme behaviour of odd order diagonally and anti-diagonally symmetric ASMs (DASASMs). We explain this briefly below.

Consider the DASASM of order $9$ in Subsection \ref{subsec:vhs-1}. We notice that this DASASM is determined by the entries in red. This portion of the matrix is called a \textit{fundamental triangle}. In general a DASASM of order $2n+1$ with entries $a_{i,j}$ $(1\leq i,j\leq 2n+1)$ is determined by the fundamental triangle $\{(i,j)|1\leq i\leq n+1, i\leq j \leq 2n+2-i\}$. Any DASASM of order $2n+1$ with $2n$ entries equal to $0$ along the portions of the diagonals that lie in this fundamental triangle is called an OOSASM of order $2n+1$. The central entry of such a matrix is always $(-1)^n$. The matrix in Subsection \ref{subsec:vhs-1} is an example of an odd order OOSASM. The aim of this section is to give generating functions for the refined enumeration of OOSASMs with respect to the position of the unique $1$ in the first row of such matrices.

\begin{figure}[!htb]
\centering
\includegraphics[scale=1.2]{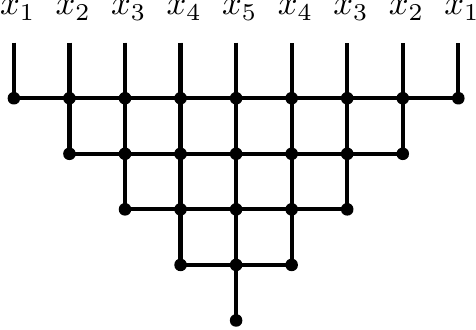}
\caption{The grid corresponding to DASASMs.}
\label{fig:grid-oos}
\end{figure}

The correspondence between ASMs and the six-vertex model extends to this case as well via the grid in Figure \ref{fig:grid-oos}. The correspondence between the degree $4$ vertices of the grid and the entries in the fundamental triangle is same as in Figure \ref{weights}. The vertices of degree $1$, namely $\raisebox{0.0005cm}{\scalebox{0.5}{\includegraphics{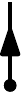}}}$ and $\raisebox{0.0005cm}{\scalebox{0.5}{\includegraphics{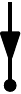}}}$ corresponds to entries $1$ and $-1$ respectively. Both of them carry a weight of $1$. The remaining vertices of degree $2$ are also fixed in our case as the boundary entries are always a $0$. So we take their weights to be $1$ as well, as they do not make any difference in enumeration results. We also assume $n\geq 1$ unless otherwise mentioned. For $n=0$, there is no OOSASM of order $1$.

Ayyer, Behrend and the first author proved the following theorem for the partition function of OOSASMs $Z_{OO}(n;x_1, x_2, \ldots, x_{n+1})$ \cite[Theorem 7.1]{AyyerBehrendFischer} (up to some normalization factor).
\begin{theo}
The OOSASM partition function of order $n$ is given by
\[
Z_{OO}(n;x_1, x_2, \ldots, x_{n+1})=Z_O\left(\left\lceil \frac{n}{2}\right\rceil;x_1, x_2, \ldots, x_{2\left\lceil \frac{n}{2} \right\rceil}\right)Q\left(\left\lceil \frac{n+1}{2}\right\rceil;x_1, x_2, \ldots, x_{2\left\lceil \frac{n+1}{2}\right \rceil-1}\right),
\]
where $Z_O(m;x_1, x_2, \ldots, x_{2m})$ is given by equation \eqref{partition-zo} and
\begin{multline*}
    Q(m;x_1, x_2, \ldots, x_{2m-1})=\sigma(q^2)^{-(m-1)(2m-1)}\prod_{1\leq i<j \leq 2m}\frac{\sigma(qx_ix_j)\sigma(q\bar x_i \bar x_j)}{\sigma(x_i \bar x_j)}\\
   \times
\underset{1\leq i<j\leq 2m}{\pf}
\left(
\begin{cases} \frac{\sigma(x_i \bar x_j)}{\sigma(qx_ix_j)} + \frac{ \sigma(x_i \bar x_j)}{\sigma(q \bar x_i \bar x_j)}, & j<2m \\
1, & j=2m \end{cases} \right).
\end{multline*}
\end{theo}

In the following, we will specialize
\[
(x_1, x_2, \ldots, x_{n+1})=(x, 1, 1, \ldots, 1) \quad \text{as well as} \quad  q+\bar q=1
\]
in the partition function. We will now explore how this specialization can be expressed in terms of $\aoo(2n+1,i)$, the number of OOSASMs of order $2n+1$ where the unique $1$ in the first row is at the $i$-th column.

We notice that there is a unique occurrence of a $\one$ in the first row, say at position $i$. This forces the other degree $4$ vertices to its left to be of type $\ne$, and to its right to be of type $\nw$. The left boundary vertex is forced to be $\raisebox{0.0005cm}{\scalebox{0.5}{\includegraphics{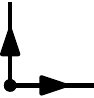}}}$ and the right boundary is forced to be $\raisebox{0.0005cm}{\scalebox{0.5}{\includegraphics{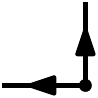}}}$. This gives us
\begin{equation}\label{oos-1}
    Z_{OO}(n;x,\underbrace{1,\ldots,1}_{n}) = \sum_{i=1}^{2n+1} \aoo(2n+1,i) \left(\frac{\sigma(qx)}{\sigma(q^2)}  \right)^{i-2} \left(\frac{\sigma(q \bar x)}{\sigma (q^2)}\right)^{2n-i}.
\end{equation}

We now perform the change of variable $$z=\frac{\sigma(q\bar x)}{\sigma(qx)}$$ and use $q+\bar q=1$ to eliminate $x$ from equation \eqref{oos-1} to get
\begin{equation}\label{oos-zo}
        Z_{OO}(n;x,\underbrace{1,\ldots,1}_{n}) = \frac{(-1)^{n-1}z^{2n}}{(z-1-z^2)^{n-1}}\sum_{i=1}^{2n+1} \aoo(2n+1,i)z^{-i}.
\end{equation}

Ayyer, Behrend and the first author \cite[Theorem 7.2]{AyyerBehrendFischer} also showed that
\begin{multline}\label{oos-1-sp}
    Z_{OO}(2n-1;x_1,x_2, \ldots , x_{2n})=3^{-(n-1)(2n-1)}Sp_{4n}(n-1,n-1,n-2,n-2,\ldots, 0,0;x_1^2, x_2^2, \ldots, x_{2n}^2)\\
    \times Sp_{4n-2}(n-1,n-2,n-2,n-3, n-3,\ldots, 1,1;x_1^2, x_2^2, \ldots, x_{2n-1}^2)
\end{multline}
and
\begin{multline}\label{oos-2-sp}
    Z_{OO}(2n;x_1,x_2, \ldots , x_{2n+1})=3^{-n(2n-1)}Sp_{4n}(n-1,n-1,n-2,n-2,\ldots, 0,0;x_1^2, x_2^2, \ldots, x_{2n}^2)\\
    \times Sp_{4n+2}(n,n-1,n-1,n-2,n-2,\ldots, 0,0;x_1^2, x_2^2, \ldots, x_{2n+1}^2).
\end{multline}

By comparing equations \eqref{eq:sp-ao}, \eqref{sp-4n+2}, \eqref{oos-zo}, \eqref{oos-1-sp} and \eqref{oos-2-sp} we get the following pairs of equations
\begin{multline}\label{oos-1-ff}
    (-1)^{n-1}3^{-n^2+2n-1}(z-1-z^2)^{n-1}\left(\sum_{i=2}^{2n}\ao(2n,i)z^{-i}\right)\left(\sum_{1\leq i\leq j\leq n}Q_{n-1,i}x^{2i-4j+2}\right)\\
    =\sum_{i=1}^{4n-1}\aoo(4n-1,i)z^{2n-2-i},
\end{multline}
where $n>1$; and
\begin{multline}\label{oos-3-ff}
    (-1)^{n}3^{-n^2}(z-1-z^2)^{n}\left(\sum_{i=2}^{2n}\ao(2n,i)z^{-i}\right)\left(\sum_{1\leq i\leq j\leq n+1}Q_{n,i}x^{2i-4j+2}\right)\\
    =\sum_{i=1}^{4n+1}\aoo(4n+1,i)z^{2n-i}.
\end{multline}

From equation \eqref{oos-1-ff} we get the following result.
\begin{theo}\label{thm:oos-1}
Let $\aoo(4n-1,i)$ denote the number of diagonally and off-diagonally symmetric ASMs of order $4n-1$ with the unique $1$ of the first row in the $i$-th column. Then for all $n>1$ the following is satisfied
\begin{multline*}
3^{-n^2+2n-1}\left(\sum_{i=2}^{2n}\ao(2n,i)z^{-i}\right)\left(\sum_{1\leq i\leq j\leq n}Q_{n-1,i}(zq-1)^{n+i-2j}(q-z)^{n-i+2j-2} (-q)^{-n+1}\right)\\
    =\sum_{i=1}^{4n-1}\aoo(4n-1,i)z^{2n-2-i},
\end{multline*}
where every quantity appearing in the left-hand side is explicitly known.
\end{theo}
Further, from equation \eqref{oos-3-ff}, we get the following result.
\begin{theo}\label{thm:oos-3}
Let $\aoo(4n+1,i)$ denote the number of diagonally and off-diagonally symmetric ASMs of order $4n+1$ with the unique $1$ of the first row in the $i$-th column. Then, for all $n\geq 1$ the following is satisfied
\begin{multline*}
3^{-n^2}\left(\sum_{i=2}^{2n}\ao(2n,i)z^{-i}\right)\left(\sum_{1\leq i\leq j\leq n+1}Q_{n,i}(zq-1)^{n+i-2j+1}(q-z)^{n-i+2j-1}(-q)^{-n}\right)\\
    =\sum_{i=1}^{4n+1}\aoo(4n+1,i)z^{2n-i},
\end{multline*}
where every quantity appearing in the left-hand side is explicitly known.
\end{theo}
\begin{rem}
From Theorems \ref{thm:vh-1} and \ref{thm:oos-1} we get the following
\begin{equation}\label{vhs-oos}
    \avh(4n+1,i)=\aoo(4n-1,i)+\aoo(4n-1,i-1),
\end{equation}
and, from Theorems \ref{thm:vh} and \ref{thm:oos-3} we get the following
\begin{equation}\label{vhs-oos-2}
    \avh(4n+3,i)=\aoo(4n+1,i)+\aoo(4n+1,i-1).
\end{equation}
In the above we assume $\aoo(2n+1,-1)=0$. These are similar to the relationship between refined enumeration of VSASMs and OSASMs (cf. equation \eqref{cni}).
\end{rem}

\begin{rem}
Analogous to equations \eqref{schur-vhs-3} and \eqref{schur-vh-1} we shall get the following pairs of equations
\begin{multline*}
    (-1)^{n^2}3^{-2n^2-n}(z-1-z^2)^{2n}
 s(2n,2n-1,2n-1,\ldots,1,1,0,0;x^2,\bar x^2,\ldots, 1)
\\=\sum_{i=1}^{4n-1}\aoo(4n-1,i)z^{4n-i}
\end{multline*}
and
\begin{multline*}
-3^{-2n^2-3n-1}(z-1-z^2)^{2n+1}s(2n+1,2n,2n,\ldots,1,1,0,0;x^2,\bar x^2,1, \ldots,1)\\
=\sum_{i=1}^{4n+1} \aoo(4n+1,i)z^{4n+2-i}.
\end{multline*}
\end{rem}

\section{Vertically and Off-Diagonally Symmetric ASMs}\label{sec:vos}

An ASM which is vertically symmetric as well as off-diagonally symmetric with a null diagonal except for the central entry is called a vertically and off-diagonally symmetric ASM (VOSASM). These matrices occur for odd orders $8n+1$ and $8n+3$. This class was first considered by Okada \cite{OkadaCharacters}, who proved enumeration formulas for them. VOSASMs are also OOSASMs of odd order as described in the previous section\footnote{The off-antidiagonal symmetry follows from the vertical symmetry, which in turn makes them vertically and horizontally symmetric. In fact VOSASMs are special cases of totally symmetric ASMs (TSASMs). No product formula for TSASMs is currently known or conjectured.}. An example of such a VOSASM is the matrix from Section \ref{subsec:vhs-1}. Clearly half of the triangular array of numbers in red is sufficient to construct the whole matrix in our example. Since VOSASMs have vertical symmetry so we can ask for their refined enumeration with respect to the position of the first occurrence of a $1$ in the second row. Let these numbers be denoted by $\avos(n,i)$ for order $n$ VOSASMs. The aim of this section is to give generating functions for these numbers.

\begin{figure}[!htb]
\centering
\includegraphics[scale=.7]{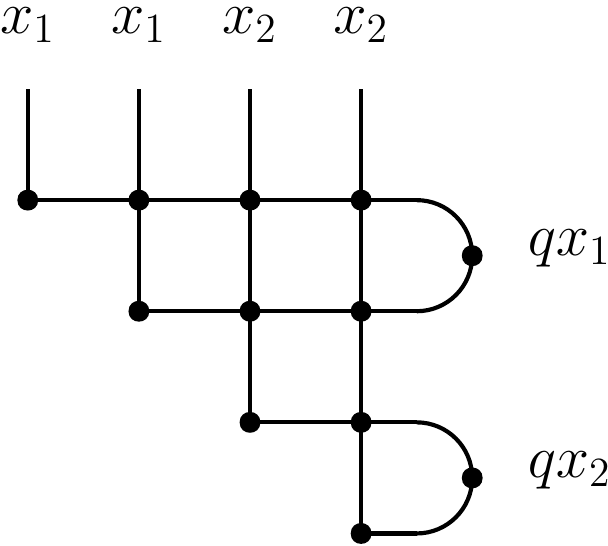}
\caption{The grid corresponding to VOSASMs.}
\label{fig:grid-vos}
\end{figure}

The correspondence between VOSASMs and the six-vertex model is via the grid shown in Figure \ref{fig:grid-vos}, which is now a combination of U-turns and the triangular grid from last section. This grid was first considered by Kuperberg \cite{KuperbergRoof}. The weights of the degree $4$ vertices remain the same as in Figure \ref{weights}, the U-turns have the weights described in Figure \ref{Uweights} and the degree $2$ vertices will have weight $1$ as they do not make any difference in our results.

We shall use the following formula for the partition function that was derived by Kuperberg \cite[Theorem 10]{KuperbergRoof} (up to some normalization factor).
\begin{theo}
\begin{multline*}
    Z_{UO}(n;x_1, x_2, \ldots, x_{2n})=\\
    \frac{\sigma(q^2)^{3n-8n^2}\sigma(q)^{2n}\prod_{i\leq 2n}\sigma(q^2x_i^2)}{\prod_{i<j\leq 2n}\sigma(\bar x_i x_j)^2\prod_{i\leq j \leq 2n}\sigma(x_ix_j)^2}\prod_{i<j\leq 2n}\sigma(q\bar x_i x_j)^2\sigma(qx_i\bar x_j)^2\sigma(qx_ix_j)^2\sigma(q\bar x_i \bar x_j)^2\\
    \times \pf_{1 \le i < j \le 2n}\left(\sigma(\bar x_i x_j)\sigma(x_ix_j)\left(\frac{1}{\sigma(qx_ix_j)\sigma(q\bar x_i \bar x_j)}-\frac{1}{\sigma(q\bar x_i x_j)\sigma(qx_i\bar x_j)}\right)\right)\\
\times \pf_{1 \le i,j \le 2n}\left(\sigma(\bar x_i x_j)\sigma(x_ix_j)\left(\frac{\sigma(bx_i)\sigma(bx_j)}{\sigma(qx_ix_j)}-\frac{\sigma(bx_i)\sigma(b\bar x_j)}{\sigma(qx_i\bar x_j)}\right.\right. \\
\left.\left.-\frac{\sigma(b\bar x_i)\sigma(bx_j)}{\sigma(q\bar x_i x_j)}+\frac{\sigma(b\bar x_i)\sigma(b\bar x_j)}{\sigma(q\bar x_i \bar x_j)}\right)\right)
\end{multline*}
\end{theo}

We shall again specialize
\[
(x_1, x_2, \ldots, x_{2n})=(x, 1, \ldots, 1) \quad \text{as well as} \quad q+ \bar q=1,
\]
in the following. But first, we notice that due to the imposed symmetries we have
\[
\avos(n,1)=\avos(n,2)=0 \quad \text{for all $n$}.
\]
There will be two cases depending on whether the VOSASM is of order $8n+1$ or $8n+3$. For order $8n+1$ we further specialize
\begin{equation}\label{special-vos-1}
   q+\bar q=1 \quad \text{as well as} \quad  b=q,
\end{equation}
and also assume $n>1$ for this case. When $n=1$, the only such VOSASM of order $9$ is the one in Section \ref{subsec:vhs-1}, and for $n=0$ no VOSASM exist. For order $8n+3$ we specialize
\begin{equation}\label{special-vos-3}
   q+\bar q=1 \quad \text{as well as} \quad  b=\bar q,
\end{equation}
and assume $n\geq 1$. When $n=0$ in this case, then the only such VOSASM is
\[
\begin{pmatrix}
 0&1&0\\
 1&-1&1\\
 0&1&0
\end{pmatrix}.
\]
The calculations for these cases are similar to the cases for VSASMs and order $4n+3$ VHSASMs. Analogous to equations \eqref{vs-eq-1} and \eqref{vsh-eq-1} we shall get the following pairs of equations.

\begin{multline*}
    Z_{UO}(n;x,\underbrace{1,\ldots,1}_{2n-1}) =  \sum_{i=3}^{4n} \avos(8n+1,i) \left( \frac{\sigma(q \bar x)}{\sigma(q x)} \right)^{i} 
\left( \frac{\sigma(q x)}{\sigma(q^2)} \right)^{8n-2} \left( \frac{\sigma(q \bar x)}{\sigma(q^2)} \right)^{-3} \sigma(q^2x) \sigma( q^2)^{2n-1}\\
+ \sum_{j=3}^{4n}  \avos(8n+1,j) \left( \frac{\sigma(q \bar x)}{\sigma(q x)} \right)^{j+2} 
\sum_{i=j+1}^{4n} \left( \frac{\sigma(q x)}{\sigma(q^2)} \right)^{2i-2} 
\left( \frac{\sigma(q \bar x)}{\sigma(q^2)} \right)^{8n-2i-6} \sigma( \bar x) \sigma( q^2)^{2n-1} \\
+ \sum_{i=3}^{4n}  \avos(8n+1,i)\left( \frac{\sigma(q x)}{\sigma(q \bar x)} \right)^{i} 
\left( \frac{\sigma(q \bar x)}{\sigma(q^2)} \right)^{8n-3} \left(\frac{\sigma(qx)}{\sigma(q^2)}\right)^{-2}\sigma( \bar x) \sigma( q^2)^{2n-1},
\end{multline*}
and
\begin{multline*}
    Z_{UO}(n;x,\underbrace{1,\ldots,1}_{2n-1}) = \sum_{i=2}^{4n} \evos(8n+3,i) \left( \frac{\sigma(q \bar x)}{\sigma(q x)} \right)^{i} 
\left( \frac{\sigma(q x)}{\sigma(q^2)} \right)^{8n-2} \left( \frac{\sigma(q \bar x)}{\sigma(q^2)} \right)^{-3} \sigma(x) \sigma(\bar q^2)^{2n-1} \\
+ \sum_{j=2}^{4n} \evos(8n+3,j) \left( \frac{\sigma(q \bar x)}{\sigma(q x)} \right)^{j+2} 
\sum_{i=j+1}^{4n} \left( \frac{\sigma(q x)}{\sigma(q^2)} \right)^{2i-2} 
\left( \frac{\sigma(q \bar x)}{\sigma(q^2)} \right)^{8n-2i-6} \sigma(\bar q^2 \bar x) \sigma(\bar q^2)^{2n-1} \\
+ \sum_{i=2}^{4n}\evos(8n+3,i) \left( \frac{\sigma(q x)}{\sigma(q \bar x)} \right)^{i} 
\left( \frac{\sigma(q \bar x)}{\sigma(q^2)} \right)^{8n-3} \left(\frac{\sigma(qx)}{\sigma(q^2)}\right)^{-2}\sigma(\bar q^2 \bar x) \sigma(\bar q^2)^{2n-1},
\end{multline*}
where
\[
\evos(8n+3,i)=\avos(8n+3,i+1)-\avos(8n+3,i).
\]

We substitute $z=\dfrac{\sigma(q\bar x)}{\sigma(qx)}$ and eliminate $x$ from the above equations to get the following pairs of equations
\begin{equation}\label{vos-1-st}
   -\sigma(q^2)^{6n-4}\sigma(q\bar x)^{-8n+4}\frac{1+z}{1-2z}Z_{UO}(n;x,\underbrace{1,\ldots,1}_{2n-1})=\sum_{i=3}^{4n}\avos(8n+1,i)(z^{i-8n+1}+z^{-i+2}),
\end{equation}
provided \eqref{special-vos-1} holds; and
\begin{equation}\label{vos-3-st}
-\sigma(q^2)^{6n-4}\sigma(q\bar x)^{-8n+4}\frac{1-z^2}{1-2z}Z_{UO}(n;x,\underbrace{1,\ldots,1}_{2n-1})=\sum_{i=2}^{4n}\evos(8n+3,i)(z^{i-8n+1}-z^{-i+3}),
\end{equation}
provided \eqref{special-vos-3} holds.

Okada \cite[Theorem 2.5]{OkadaCharacters} proved that
\begin{multline}\label{okada-vos-1}
    \prod_{i=1}^{2n}(x_i+\bar x_i)\prod_{i=1}^{2n}\sigma(q^2x_i^2)^{-1}Z_{UO}(n;x_1, x_2, \ldots, x_{2n})\\=3^{-4n^2+3n}
   ( Sp_{4n}(n-1,n-1,n-2,n-2,\ldots,0,0;x_1^2,\ldots,x_{2n}^2))^3\\ 
   \times O_{4n}\left(n+\frac{1}{2},n-\frac{1}{2},n-\frac{1}{2},n-\frac{3}{2},n-\frac{3}{2},\ldots,\frac{3}{2},\frac{3}{2},\frac{1}{2};x_1^2,\ldots,x_{2n}^2 \right)
\end{multline}
provided \eqref{special-vos-1} holds; and
\begin{multline}\label{okada-vos-3}
    \prod_{i=1}^{2n}\sigma(q^2x_i^2)^{-1}Z_{UO}(n;x_1, x_2, \ldots, x_{2n})\\=3^{-4n^2+3n}
   ( Sp_{4n}(n-1,n-1,n-2,n-2,\ldots,0,0;x_1^2,\ldots,x_{2n}^2))^3\\ \times Sp_{4n+2}(n,n-1,n-1,n-2,n-2,\ldots,0,0;x_1^2,\ldots,x_{2n}^2,1)
\end{multline}
provided \eqref{special-vos-3} holds. Transforming the even orthogonal group character in equation \eqref{okada-vos-1} for our special case, into a symplectic group character using a result of Ayyer and Behrend \cite[Proposition 5, and then use equation (7)]{AyyerBehrend}, we get
\begin{multline}\label{okada-vos-2}
 Z_{UO}(n;x,\underbrace{1,\ldots,1}_{2n-1}) =3^{-4n^2+5n-2}\sigma(q^2)^{2n-1}\sigma(q^2x^2)(x^2+1+\bar x^2)\\
\times    (Sp_{4n}(n-1,n-1,n-2,n-2,\ldots,0,0;x^2,1,\ldots,1))^3 \\
\times Sp_{4n-2}(n-1,n-2,n-2,\ldots,0,0;x^2,1,\ldots,1).\\ 
\end{multline}
provided \eqref{special-vos-1} holds. 

From equations \eqref{eq:sp-ao}, \eqref{sp-4n+2}, \eqref{vos-1-st} and \eqref{okada-vos-2}, after some simplifications we get
\begin{multline}\label{vos-1-final}
(-1)^{n-1}3^{-n^2+2n-1}(1+z)(z-1-z^2)^{n-1}\left(\sum_{i=2}^{2n}\ao(2n,i)z^{-i}\right)^3 \left( \sum_{1 \le j \le i \le n} Q_{n-1,i} x^{2i-4j+2}\right) \\ =
\sum_{i=3}^{4n}\avos(8n+1,i)(z^{i-6n-4}+z^{2n-i-3}),
\end{multline}
provided \eqref{special-vos-1} holds and where $Q_{n-1,i}$ is given by \eqref{value-qni}. On the other hand, from equations \eqref{eq:sp-ao}, \eqref{sp-4n+2}, \eqref{vos-3-st} and \eqref{okada-vos-3} we get
\begin{multline}\label{vos-3-final}
(-1)^n3^{-n^2}(1-z^2)z^{-2n+4}(z-1-z^2)^{n}
\left(\sum_{i=2}^{2n}\ao(2n,i)z^{-i}\right)^3\left(\sum_{1\leq i\leq i\leq n+1}Q_{n,i}x^{2i-4j+2}\right)\\
=\sum_{i=2}^{4n}\evos(8n+3,i)(z^{i-8n+1}-z^{-i+3}),
\end{multline}
provided \eqref{special-vos-3} holds and where $Q_{n,i}$ is given by \eqref{value-qni}.

After some simplifications the above pairs of equations give the following theorems.
\begin{theo}\label{thm:vos-1}
Let $\avos(8n+1,i)$ denote the number of order $8n+1$ VOSASMs with the first $1$ in the second row in the $i$-th column. Then, for all $n>1$ the following is satisfied
\begin{multline*}
3^{-n^2+2n-1}(1+z)\left(\sum_{i=2}^{2n}\ao(2n,i)z^{-i}\right)^3 \left( \sum_{1 \le j \le i \le n} Q_{n-1,i} (zq-1)^{n+i-2j}(q-z)^{n-i+2j-2}(-q)^{-n+1}\right) \\ =
\sum_{i=3}^{4n}\avos(8n+1,i)(z^{i-6n-4}+z^{2n-i-3}),
\end{multline*}
where every quantities appearing on the left hand side is explicitly known.
\end{theo}
\begin{theo}\label{thm:vos-3}
Let $\avos(8n+3,i)$ denote the number of order $8n+3$ VOSASMs with the first $1$ in the second row in the $i$-th column. Then, for all $n\geq 1$ the following is satisfied
\begin{multline*}
3^{-n^2}(1-z^2)
\left(\sum_{i=2}^{2n}\ao(2n,i)z^{-i}\right)^3\left(\sum_{1\leq i\leq i\leq n+1}Q_{n,i}(zq-1)^{n+i-2j+1}(q-z)^{n-i+2j-1}(-q)^{-n}\right)\\
=\sum_{i=2}^{4n}(\avos(8n+3,i+1)-\avos(8n+3,i))(z^{i-6n-3}-z^{2n-i-1}),
\end{multline*}
where every quantities appearing on the left hand side is explicitly known, and $\avos(8n+3,1)=\avos(8n+3,2)=0$ for all $n$.
\end{theo}

\begin{rem}
From Theorems \ref{thm:vh-1} and \ref{thm:vos-1} we get
\begin{multline}
    \label{rel-vhs-vos-1}
   \left (\sum_{i=2}^{2n}\ao(2n,i)z^{-i}\right)^2\left(\sum_{i=2}^{2n}\avh(4n+1,i)(z^{i-2n-2}+z^{2n-1-i})\right)\\
   =\sum_{i=3}^{4n}\avos(8n+1,i)(z^{i-6n-4}+z^{2n-i-3})
\end{multline}
and from Theorems \ref{thm:vh} and \ref{thm:vos-3} we get
\begin{multline}
    \label{rel-vhs-vos-3}
   \left (\sum_{i=2}^{2n}\ao(2n,i)z^{-i}\right)^2\left(\sum_{i=1}^{2n+1}  \left( \avh(4n+3,i+1) - \avh(4n+3,i) \right) \left( z^{i-2n-1} - z^{-i+2n+1} \right) \right)\\=\sum_{i=2}^{4n}(\avos(8n+3,i+1)-\avos(8n+3,i))(z^{i-6n-3}-z^{2n-i-1}).
\end{multline}
\end{rem}

\section{Quarter Turn Symmetric ASMs}\label{sec:qtsasm}

ASMs that are invariant under a $90^\circ$ rotation are called quarter turn symmetric ASMs (QTSASMs). As a first observation, we see that these ASMs cannot occur for order $4n+2$~\cite[Lemma 4]{aval-duchon}, consider the QTSASM of order $2n$ where the entries are given by $a_{i,j}$ ($1\leq i,j \leq 2n$). Then we have
\[
2n=\sum_{1\leq i,j \leq 2n}a_{i,j}=4\sum_{1\leq i,j\leq n}a_{i,j},
\]
and this implies that $2 |n$. So for the even case they occur only for order $4n$.  Order $4n$ QTSASMs were enumerated by Kuperbeg \cite{KuperbergRoof}, while Razumov and Stroganov \cite{rs-odd-quarter} enumerated them for odd order. In this section, we will give refined enumeration formulas for this class of ASMs with respect to the position of the unique $1$ in the first row. These results were conjectured by Robbins \cite{robbins-2} and are almost immediate from several results that were already proved by Kuperberg \cite{KuperbergRoof} and by Razumov and Stroganov \cite{rs-odd-quarter}, but which to the best of our knowledge have not appeared explicitly in the literature so far. We deal with the even case in Subsection \ref{sub:even} and with the odd case in Subsection \ref{sub:odd}.

Since the proofs are similar for all cases (in this section and also in Section \ref{sec:qqtsasm}), we will only derive the results in full for the even QTSASM case and indicate the steps for the other cases. Throughout this section, we will consider the case $\vec{x}=(x,1,1,\ldots, 1)$ and $q+\bar q=1$, unless otherwise mentioned.

\subsection{Even Order QTSASMs}\label{sub:even}

\begin{figure}[!htb]
\centering
\includegraphics[scale=.6]{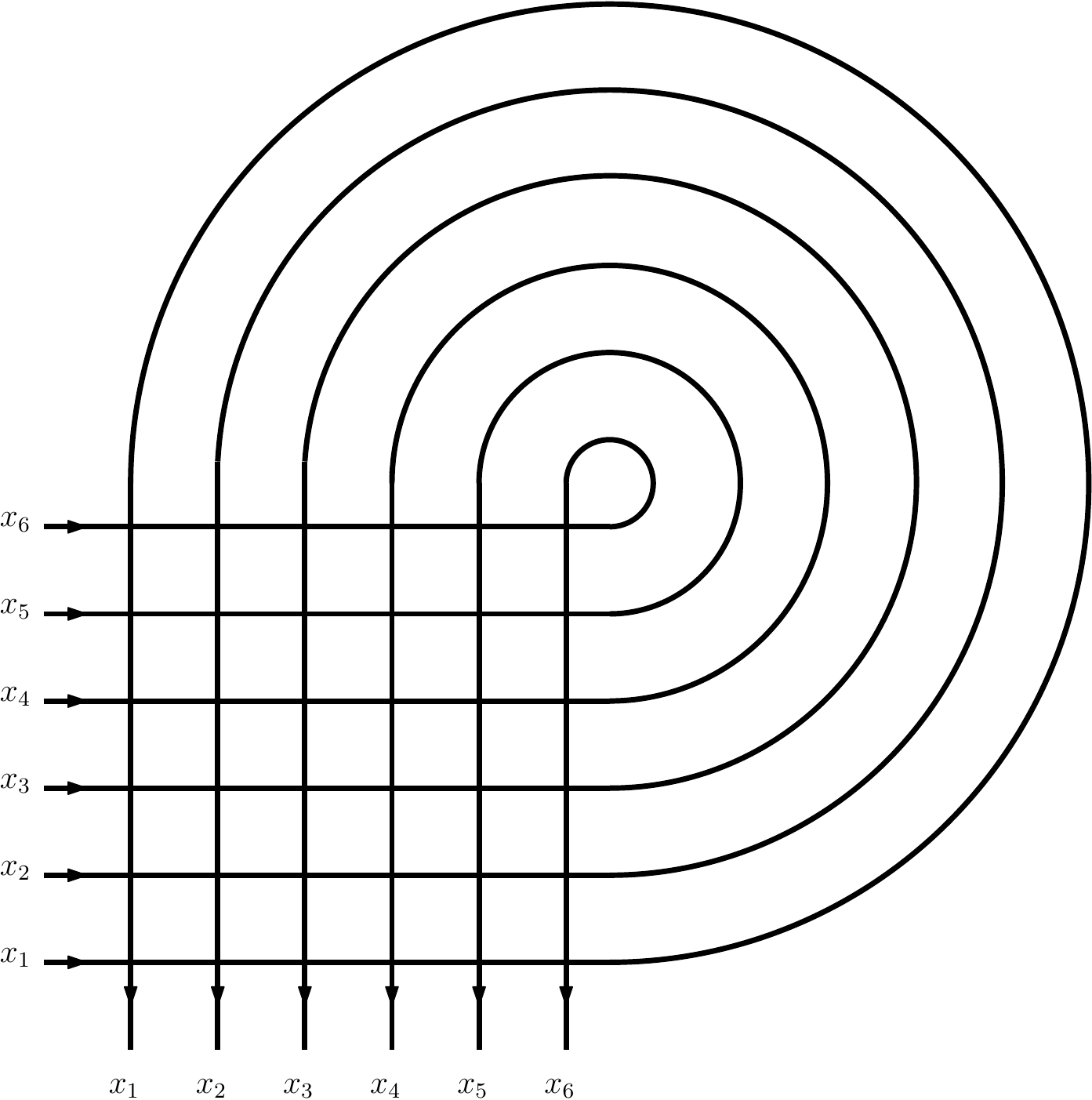}
\caption{Six vertex configuration of QTSASM of order $4n$.}
\label{fig:qt-4m}
\end{figure}

Kuperberg \cite{KuperbergRoof} enumerated this class of QTSASMs using the six vertex model; since a quarter of such an ASM is enough to determine the whole ASM, we can see that the six vertex model corresponding to even order QTSASMs will be the one shown in Figure \ref{fig:qt-4m}, where there are $2n$ many spectral parameters associated with it if the QTSASM is of order $4n$, which are denoted by $x_i$'s. The vertex weights and bijection with the normal ASMs which was shown in Section \ref{asm} carry forward to this class as well, with one difference: the arrows of the configuration will change sign when they move through the circular turns. For instance, consider the simplest possible QTSASM of even order 
\[
\begin{pmatrix}
 0&1&0&0\\
 0&0&0&1\\
 1&0&0&0\\
 0&0&1&0
\end{pmatrix}.
\]
The six-vertex configuration of this matrix is now given by Figure \ref{fig:qt-ex-eps-converted-to.pdf}.
\begin{figure}[!htb]
\centering
\includegraphics[scale=.6]{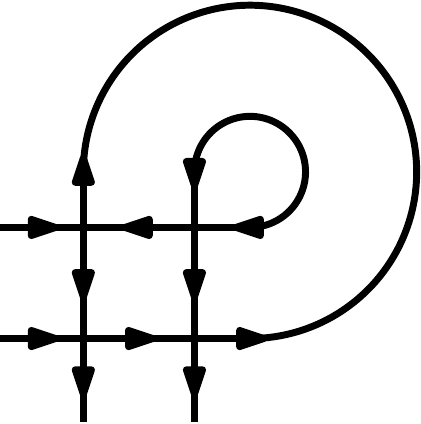}
\caption{Six vertex configuration of an even order QTSASM.}
\label{fig:qt-ex}
\end{figure}

Let $\aqt(n,i)$ be the total number of QTSASMs of order $n$ such that the unique $1$ in the last row is at the $i$-th position. By symmetry of the ASMs, this is also the number of QTSASMs of order $n$ such that the unique $1$ in the first row is at the $i$-th position. In the sequel we will use this description instead of the last row, without further commentary. The weight of the last row of such a QTSASM of order $4n$ would be \[\left(\frac{\sigma(q\overline{x})}{\sigma(q^2)}\right)^{i-2}\left(\frac{\sigma(qx)}{\sigma(q^2)}\right)^{4n-i-1}.\] (The first entry is always a $0$ for QTSASMs, and the weight of that entry is $\frac{\sigma(qx\bar x)}{\sigma(q^2)}=1$ when $q+\bar q=1$.) If $Z_{Q}(4n;x_1, x_2, \ldots, x_{2n})$ denotes the partition function of QTSASMs of order $4n$, then we have
\begin{equation}\label{eq-qt-4n}
    \sum_{i=2}^{4n-1}\aqt(4n,i)z^{i-2}=\frac{Z_{Q}(4n;x,1,1,\ldots,1)}{\sigma(q^2)^{3-4n}\sigma(qx)^{4n-3}},
\end{equation}
where $z=\dfrac{\sigma(q\overline{x})}{\sigma(qx)}$ and from Kuperberg \cite[Theorem 10]{KuperbergRoof} we have (up to some normalization factor)
\begin{multline*}
    Z_{Q}(4n; x_1, x_2, \ldots, x_{2n})=\sigma(q^2)^{4n-4n^2}\frac{\prod_{i,j=1}^{2n}\alpha(\bar x_i x_j)^2}{\prod_{i,j=1}^{2n}\sigma(\bar x_i x_j)^2}\\\pf_{1 \le i < j \le 2n}\left(\frac{\sigma(\bar x_i x_j)}{\alpha(\bar x_i x_j)}\right)\pf_{1 \le i < j \le 2n}\left(\frac{\sigma(\bar x_i^2 x_j^2)}{\alpha(\bar x_i x_j)}\right).
\end{multline*}

Further we have the following analogous equations from Razumov and Stroganov~\cite[Equation (38)]{rs-odd-half},
\begin{equation}\label{eq-asm-ref}
    \sum_{i=1}^{n}\an(n,i)z^{i-1}=\frac{Z(n;x,1,1,\ldots,1)}{\sigma(q^2)^{1-n}\sigma(qx)^{n-1}}
\end{equation}
and 
\begin{equation}\label{eq-ht-ref}
        \sum_{i=1}^{2n}\aht(2n,i)z^{i-1}=\frac{Z_{H}(2n;x,1,1,\ldots,1)}{\sigma(q^2)^{1-2n}\sigma(qx)^{2n-1}}
\end{equation}
where $\an(n,i)$ denotes the number of ASMs of order $n$ with the unique $1$ in the first row in the $i$-th column, given by Zeilberger \cite{zeil}
\[
\an(n,i)=\binom{n+i-2}{n-1}\frac{(2n-i-1)!}{(n-i)!}\prod_{j=0}^{n-2}\frac{(3j+1)!}{(n+j)!}
\]
and $\aht(n,i)$ denotes the number of order $n$ half-turn symmetric ASMs (HTSASMs)\footnote{ASMs which are invariant under a $180^\circ$ rotation are called half-turn symmetric ASMs.} with the unique $1$ in the first row in the $i$-th column, given by Stroganov \cite{stroganov-izergin}
\begin{multline*}
    \aht(2n,i)=\frac{(2n-1)!^2}{(n-1)!^2(3n-3)!(3n-1)!}\prod_{j=0}^{n-1}\frac{(3j+2)(3j+1)!^2}{(3j+1)(n+j)!^2}\\
    \times \sum_{j=1}^i\left(\frac{(n^2-nj+(j-1)^2(n+j-3)!)(2n-j-1)!(n+i-j-1)!(2n-i+j-2)!}{(j-1)!(n-j+1)!(i-j)!(n-i+j-1)!}\right).
\end{multline*}
Also $$Z(n;x_1, x_2, \ldots, x_n, y_1, y_2, \ldots, y_n)$$
\noindent is the partition function of order $n$ ASMs and 
$$Z_{H}(2n;x_1, x_2, \ldots, x_n, y_1, y_2, \ldots, y_n)$$ is the partition function of order $2n$ HTSASMs, which are given by Kuperberg~\cite[Theorem 10]{KuperbergRoof} as follows (up to some normalization factor)
\[
Z(n;x_1, \ldots, x_n, y_1, \ldots, y_n)=\frac{\sigma(q^2)^{n-n^2}\prod_{i,j}\alpha(x_i\bar x_j)}{\prod_{i<j}\sigma(\bar x_i x_j)\sigma(y_i \bar y_j)}\det_{1\leq i,j\leq n}\left(\frac{1}{\alpha(x_i\bar x_j)}\right),
\]
and
\begin{multline*}
    Z_{H}(2n;x_1, \ldots, x_n, y_1, \ldots, y_n)=\frac{\sigma(q^2)^{n-2n^2}\prod_{i,j}\alpha(x_i\bar x_j)^2}{\prod_{i<j}\sigma(\bar x_i x_j)^2\sigma(y_i \bar y_j)^2}\\\det_{1\leq i,j\leq n}\left(\frac{1}{\alpha(x_i\bar x_j)}\right)\det_{1\leq i,j \leq n}\left(\frac{1}{\sigma(q\bar x_i y_j)}+\frac{1}{\sigma(qx_i \bar y_j}\right).
\end{multline*}

By results of Okada \cite[Theorems 2.4 and 2.5]{OkadaCharacters}, when $q+\bar q=1$, we have
\begin{equation}\label{eq-qt-rel}
    Z_{Q}(4n;x, 1, \ldots ,1)=(Z(n;x,1,\ldots, 1))^2Z_{H}(2n;x,1,\ldots, 1).
\end{equation}
Combining equations \eqref{eq-qt-4n} to \eqref{eq-qt-rel}, we get the following result.
\begin{theo}\label{qt-ref-thm}
\[
\sum_{i=2}^{4n-1}\aqt(4n,i)z^{i-2}=\left( \sum_{i=1}^{n}\an(n,i)z^{i-1}\right )^2\left(\sum_{i=1}^{2n}\aht(2n,i)z^{i-1}\right).
\]
\end{theo}

Theorem \ref{qt-ref-thm} allows us to give a formula for the refined enumeration of even order QTSASMs by comparing the coefficients of $z$ from both sides of the equation.

\subsection{Odd Order QTSASMs}\label{sub:odd}

\begin{figure}[!htb]
\centering
\includegraphics[scale=.6]{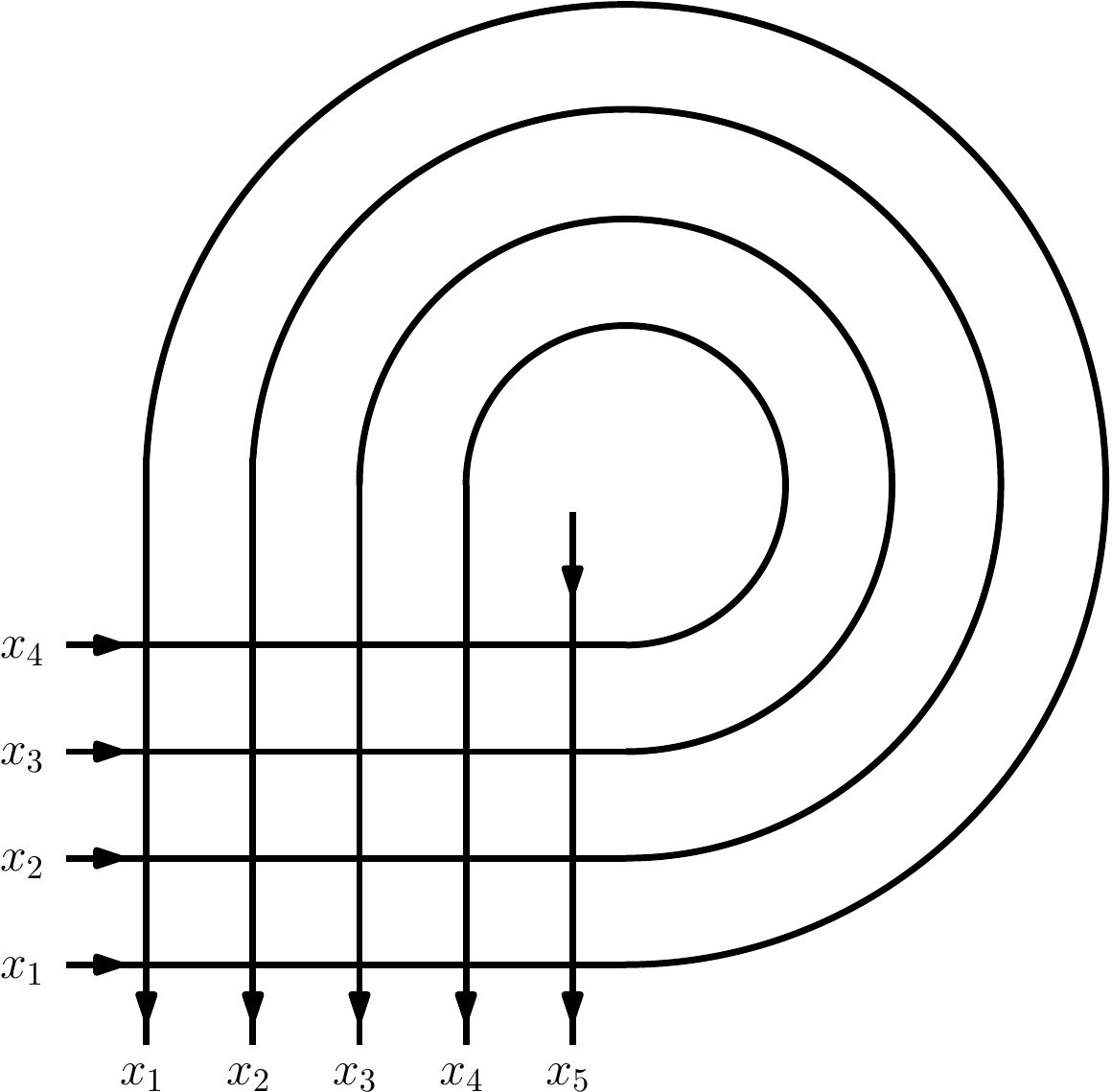}
\caption{Six vertex configuration of QTSASM of order $4n+1$.}
\label{fig:qt-4m-1}
\end{figure}

\begin{figure}[!htb]
\centering
\includegraphics[scale=.6]{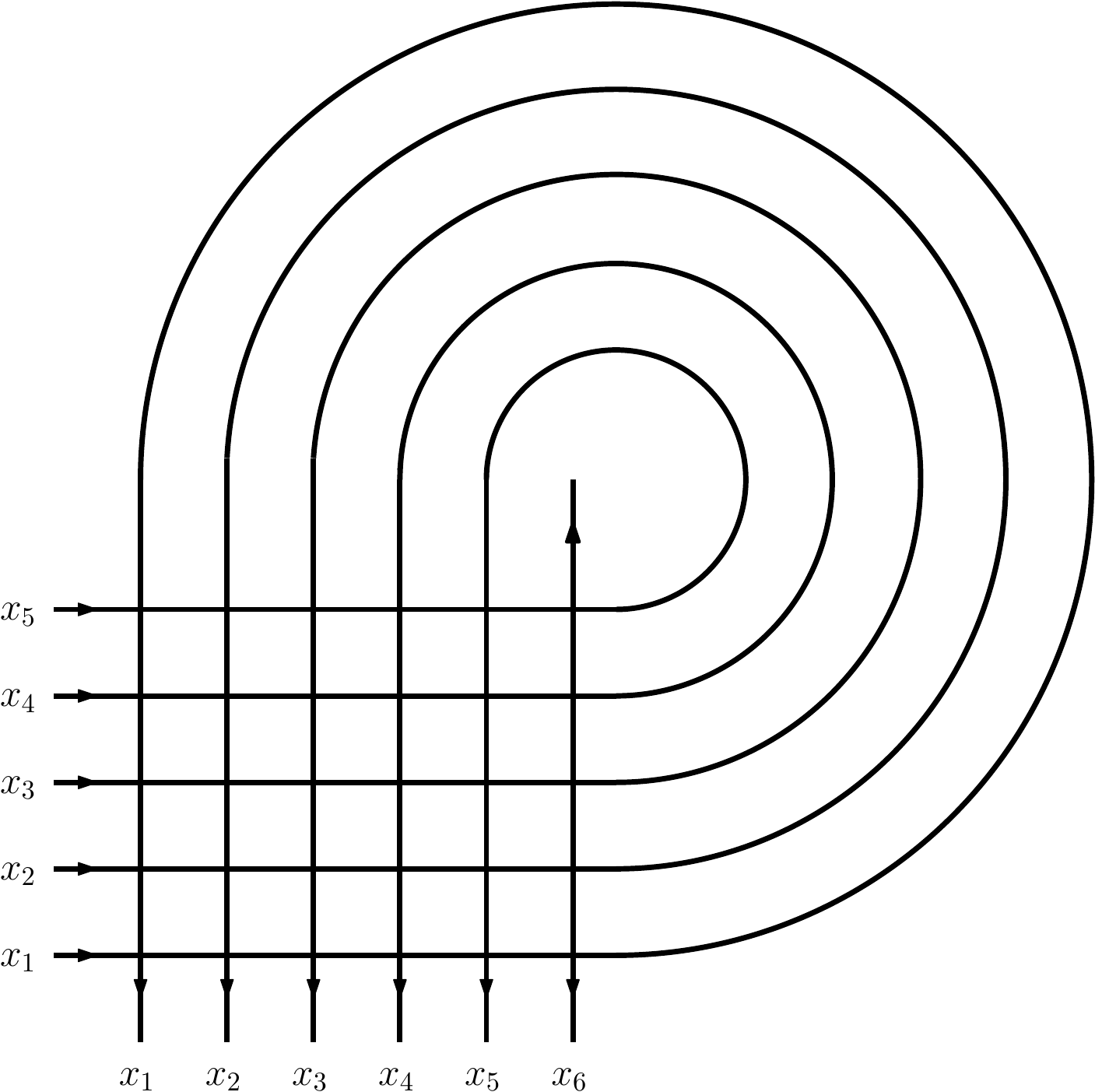}
\caption{Six vertex configuration of QTSASM of order $4n+3$.}
\label{fig:qt-4m-3}
\end{figure}

For odd order QTSASMs (say $2m+1$), we can notice that the central entry of the matrix will be $(-1)^m$. Similar to the grid described in Subsection \ref{sub:even}, we will have the configuration in Figure \ref{fig:qt-4m-1} if $m=2n$ and in Figure \ref{fig:qt-4m-3} if $m=2n+1$. These models were used by Razumov and Stroganov to enumerate odd order QTSASMs \cite{rs-odd-quarter}, and they found a  formula similar to the one found by Kuperberg for even order QTSASMs. We also notice that for an order $2m+1$ QTSASM, the grid has $m+1$ spectral parameters.

We proceed in a similar way, as in the case for even QTSASMs. Analogous to equation \eqref{eq-qt-4n}, we will have the following equations.
\begin{equation}\label{eq-odd-11-ref}
        \sum_{i=2}^{4n}\aqt(4n+1,i)z^{i-2}=\frac{Z_{Q}(4n+1;x,1,1,\ldots,1)}{\sigma(q^2)^{2-4n}\sigma(qx)^{4n-2}},
\end{equation}
and
\begin{equation}\label{eq-odd-2-ref}
        \sum_{i=2}^{4n+2}\aqt(4n+3,i)z^{i-2}=\frac{Z_{Q}(4n+3;x,1,1,\ldots,1)}{\sigma(q^2)^{-4n}\sigma(qx)^{4n}},
\end{equation}
where $z=\dfrac{\sigma(q\bar x)}{\sigma(qx)}$ and the partition functions are given by Razumov and Stroganov \cite[Equations (10) and (11)]{rs-odd-quarter} (up to some normalization factor). Further we have the following analogous equation from Razumov and Stroganov \cite{rs-odd-half},
\begin{equation}\label{eq-ht-odd-ref}
     \sum_{i=1}^{2n+1}\aht(2n+1,i;w)t^{i-1}=\frac{Z_{H}(2n+1;x,1,1,\ldots,1)}{\sigma(q^2)^{-2n}\sigma(qx)^{2n}},
\end{equation}
where $Z_{H}(2n+1; x_1, x_2, \ldots, x_{n+1},y_1, y_2, \ldots, y_{n+1})$ is the partition function of odd order HTSASMs, which was found by Razumov and Stroganov \cite[Theorem 1]{rs-odd-half} (up to some normalization factor).

Again, analogous to equation \eqref{eq-qt-rel}, we have the following equations from Razumov and Stroganov \cite{rs-odd-quarter}, when $q+\bar q=1$
\begin{equation}\label{eq-odd-qt-rs-1}
        Z_{Q}(4n+1;x,1,\ldots, 1)=(Z(n;x,1,\ldots, 1))^2Z_{H}(2n+1;x,1,\ldots,1),
\end{equation}
and
\begin{equation}\label{eq-odd-qt-rs-2}
        Z_{Q}(4n+3;x,1,\ldots, 1)=(Z(n+1;x,1,\ldots, 1))^2Z_{H}(2n+1;x,1,\ldots,1).
\end{equation}
Now, combining equations \eqref{eq-asm-ref} and \eqref{eq-odd-11-ref} to \eqref{eq-odd-qt-rs-2}, we get the following result.
\begin{theo}
\[
\sum_{i=2}^{4n}\aqt(4n+1,i)z^{i-2}=\left( \sum_{i=1}^{n}\an(n,i)z^{i-1}\right )^2\left(\sum_{i=1}^{2n+1}\aht(2n+1,i)z^{i-1}\right),
\]
and
\[
\sum_{i=2}^{4n+2}\aqt(4n+3,i)z^{i-2}= \left( \sum_{i=1}^{n+1}\an(n+1,i)z^{i-1}\right )^2\left(\sum_{i=1}^{2n+1}\aht(2n+1,i)z^{i-1}\right).
\]
\end{theo}

\section{Quasi-Quarter Turn Symmetric ASMs}\label{sec:qqtsasm}

\begin{figure}[!htb]
\centering
\includegraphics[scale=.6]{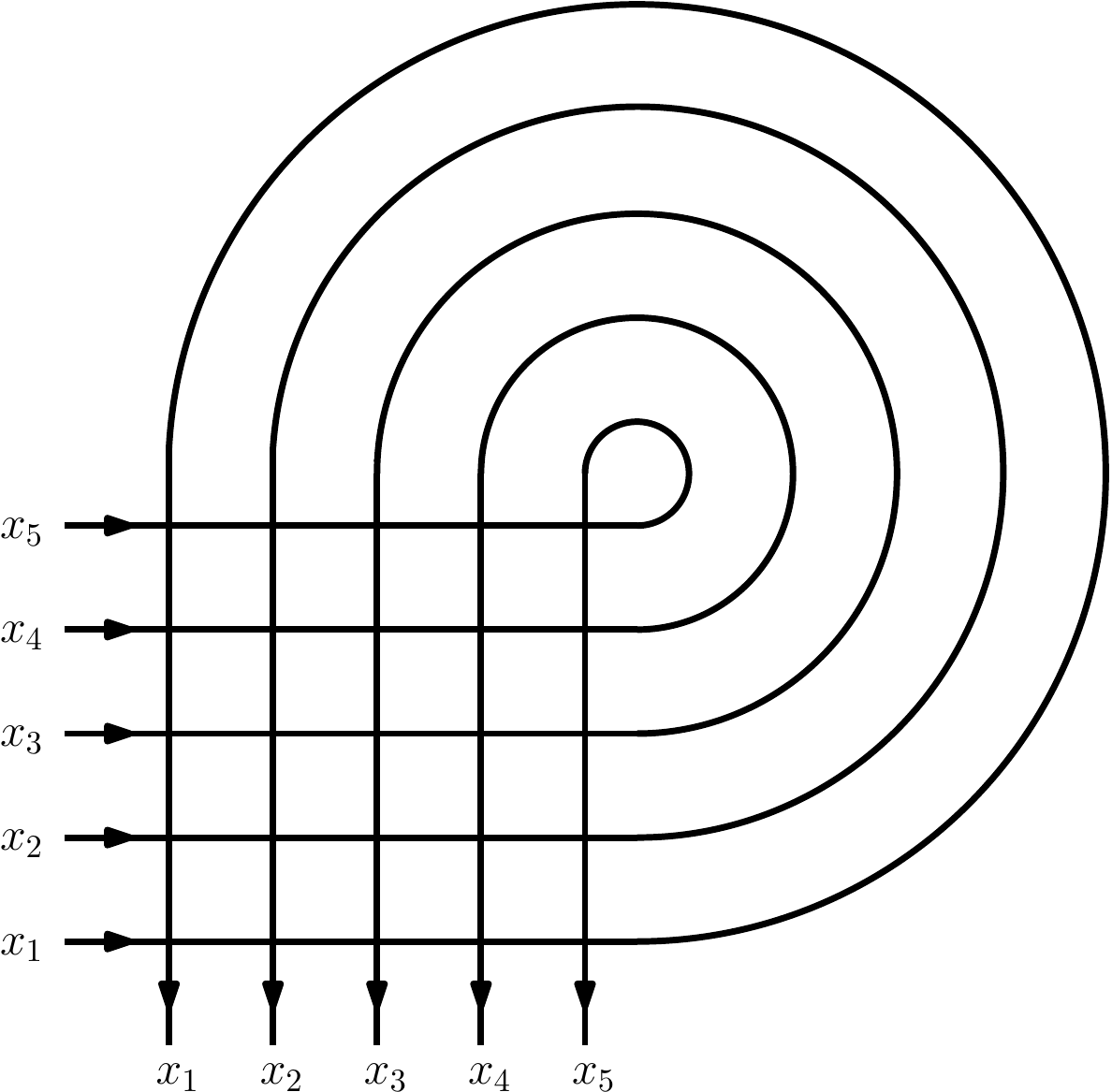}
\caption{Six vertex configuration of qQTSASM of order $4n+2$.}
\label{fig:q-qt}
\end{figure}

As pointed out in Section \ref{sec:qtsasm}, there are no even order QTSASMs of order $4n+2$. However, Duchon~\cite{duchon} introduced a new type of ASM, called quasi-QTSASMs (qQTSASMs) which follows all the conditions of an ASM and has quarter turn symmetry for all entries except the four entries in the middle, which can be either $\{1,0,0,1\}$ or $\{0,-1,-1,0\}$. Below we give examples of both these types of matrices.

\[
\begin{pmatrix}
 0&0&1&0&0&0\\
 0&1&-1&0&1&0\\
 0&0&1&0&-1&1\\
 1&-1&0&1&0&0\\
 0&1&0&-1&1&0\\
 0&0&0&1&0&0
\end{pmatrix}
\qquad  \qquad
\begin{pmatrix}
 0&0&0&1&0&0\\
 0&0&1&0&0&0\\
 1&0&0&-1&1&0\\
 0&1&-1&0&0&1\\
 0&0&0&1&0&0\\
 0&0&1&0&0&0
\end{pmatrix}
\]

Duchon \cite{duchon} stated two conjectures related to the enumeration of qQTSASMs: the first of these is on the unrestricted enumeration of this class of ASMs, which was subsequently proved by Aval and Duchon \cite{aval-duchon} by studying the six vertex configuration associated with the even order qQTSASMs (shown in Figure \ref{fig:q-qt}); the second conjecture deals with the refined enumeration of qQTSASMs with respect to the position of the unique $1$ in the first row, which we prove below.

\begin{theo}\label{qqt-ref-thm}
Let $\aqqt(n,i)$ denote the number of order $n$ qQTSASMs with the unique $1$ in the first row in the $i$-th column, then we have
\[
\sum_{i=2}^{4n+1}\aqqt(4n+2,i)z^{i-2}=\left(\sum_{i=1}^{n}\an(n,i)z^{i-1}\right)\left(\sum_{i=1}^{n+1}\an(n+1,i)z^{i-1}\right)\left(\sum_{i=1}^{2n+1}\aht(2n+1,i)z^{i-1}\right).
\]
\end{theo}

\begin{proof}
The proof is similar to the proof of Theorem \ref{qt-ref-thm}, hence we will just sketch it. Analogous to equation \eqref{eq-qt-4n} we have the following equation in this case, provided $q+\bar q=1$
\begin{equation}\label{eq-qqt}
    \sum_{i=2}^{4n+1}\aqqt(4n+2,i)z^{i-2}=\frac{Z_{qQ}(4n+2;x,1,1,\ldots,1)}{\sigma(q^2)^{1-4n}\sigma(qx)^{4n-1}},
\end{equation}
where $Z_{qQ}(4n+2;x_1, x_2, \ldots, x_{2m+1})$ is the partition function of order $4n+2$ qQTSASMs with spectral parameters $\vec{x}$, which was found by Aval and Duchon \cite{aval-duchon} (up to some normalization factor).

We now use the following result of Aval and Duchon \cite[Theorem 6]{aval-duchon}, if $q+\bar q=1$ then we have
\begin{equation}\label{ad-eq}
Z_{qQ}(4n+2;x, 1, \ldots, 1)=Z(n;x, 1, \ldots, 1)Z(n+1;x, 1, \ldots, 1)Z_{H}(2n+1;x, 1, \ldots, 1).
\end{equation}
Combining equations \eqref{eq-qqt},\eqref{eq-asm-ref},\eqref{eq-ht-ref} and \eqref{ad-eq} we obtain the result.
\end{proof}

\section{Concluding Remarks}\label{sec:rem}

In this paper, we have proved several results for refined enumeration of symmetry classes of ASMs. We summarize the results so far in the literature for singly refined enumeration of symmetry classes of ASMs below.
\begin{itemize}
    \item {\it Unrestricted ASMs.} First enumerated by Zeilberger \cite{doron}. The refined enumeration with respect to first row is also due to Zeilberger \cite{zeil}.
    \item {\it VSASMs.} Enumerated by Kuperberg \cite{KuperbergRoof}, singly refined enumeration with respect to the first column was accomplished by Razumov and Stroganov \cite{rs-refined}, while the singly refined enumeration with respect to the second row is obtained in this paper.
    \item {\it VHSASMs.} Enumerated by Okada \cite{OkadaCharacters}. Generating function for singly refined enumeration with respect to the second row provided in this paper.
    \item {\it HTSASMs.} For even order enumerated by Kuperberg \cite{KuperbergRoof} and for odd order by Razumov and Stroganov \cite{rs-odd-half}. Singly refined enumeration done by Stroganov \cite{stroganov-izergin}.
    \item {\it QTSASMs.} Enumerated by Kuperberg \cite{KuperbergRoof} for order $n\equiv 0 \pmod 4$ and by Razumov and Stroganov \cite{rs-odd-half} for order $n\equiv 1 \pmod 2$. Refined enumeration results proved in this paper.
    \item {\it Diagonally Symmetric ASMs.} No formula known or conjectured.
    \item {\it DADSASMs.} For odd order enumerated by Behrend, Fischer and Konvalinka \cite{bfk}.
    \item {\it TSASMs.} No formula known or conjectured.
\end{itemize}
As can be seen from the list, our results almost exhaust the singly refined enumeration possibilities for the symmetry classes of ASMs. Double refinement by fixing the positions of the $1$s in two rows or columns have been studied only for the case of unrestricted ASMs. Some other doubly refined enumeration possibilities still remain, for instance doubly refined enumeration of VSASMs with the positions of the first $1$ in the second row and the unique $1$ in the first column fixed. However, no conjectures for such cases are known at the moment. We refer the reader to Behrend, Fischer and Konvalinka's paper \cite[Section 1.2]{bfk} for a survey of all such known or conjectured results.

In addition to the symmetry classes, the following closely related classes have been studied in the literature (as well as this paper).
\begin{itemize}
    \item {\it qQTSASMs.} Enumerated by Aval and Duchon \cite{aval-duchon}. Refined enumeration results proved in this paper.
    \item {\it OSASMs.} For even order enumerated by Kuperberg \cite{KuperbergRoof} and refined enumeration by Razumov and Stroganov \cite{rs-refined}.
    \item {\it OOSASMs.}  For even order no formula is known or conjectured. For odd order enumerated by Ayyer, Behrend and Fischer \cite{AyyerBehrendFischer}, refined enumeration results in terms of generating function proved in this paper.
    \item {\it VOSASMs.} Enumerated by Okada \cite{OkadaCharacters} and refined enumeration results in terms of generating function proved in this paper.
\end{itemize}

We would like to close the paper by posing the following problems.
\begin{enumerate}
    \item Bijective proofs of equations \eqref{cni}, \eqref{vhp-rel}, \eqref{vhs-oos} and \eqref{vhs-oos-2} would be interesting.
    \item The number $\avc(2n+1,i)$ of vertically symmetric ASMs with the unique $1$ in the first column situated in the $i$-th row equals $\ao(2n,i)$ (this follows from the work of Razumov and Stroganov \cite{RazumovStroganov}). This gives us
    \[
    \av(2n+1,i)=\avc(2n+1,i)+\avc(2n+1,i+1).
    \]
    Again, a bijective proof of this would be of interest.
    \item Simpler enumeration formulas for $\avh(2n+1,i)$, $\aoo(2n+1,i)$ and $\avos(n,i)$ are desirable. In fact, since the formulas for all of these refinements are related (cf. relations \eqref{vhs-oos}, \eqref{vhs-oos-2}, \eqref{rel-vhs-vos-1} and \eqref{rel-vhs-vos-3}), so a simpler formula for anyone of them would yield the others.
\end{enumerate}

\appendix

\section{Evaluation of the Symplectic Character}\label{appen}

Symplectic characters can be interpreted as a generating function of weighted rhombus tilings of certain regions in the triangular lattice, which follows from the work of Proctor \cite{proctor} and Cohn, Larsen and Propp \cite{CohnLarsenPropp} (see also the recent work of Ayyer and the first author \cite{ilse-sp}). We shall use this description now to study the symplectic character that appears in \eqref{vsh-3}.

\begin{figure}
    \centering
    \includegraphics[width=0.2\textwidth]{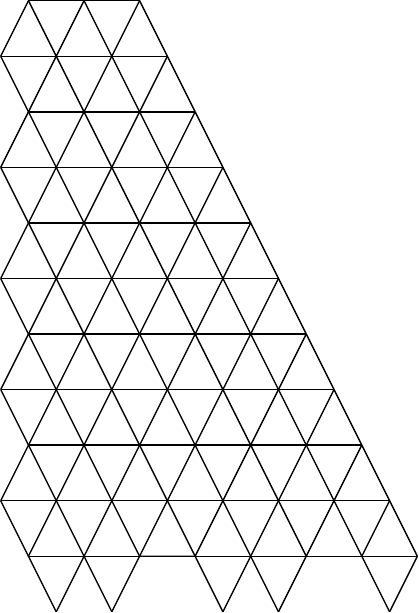}
    \caption{$\textup{QH}^{1, 2,4,5,7}_{10,2}(x_1, x_2, \ldots, x_{10})$.}
    \label{fig:ht-2}
\end{figure}

We start with a quartered hexagon in the triangular lattice, which is obtained after cutting it along its axes, see Figure~\ref{fig:ht-2}. Let the top boundary of this quartered hexagon be of length $\ell$ and the right boundary of length $2n$. The left boundary consists of a zig-zag line containing $n$ up-pointing and $n$ down-pointing triangles, and the bottom boundary has protruded down-pointing triangles in positions $p_1<p_2<\cdots <p_n$, numbered from left to right, starting with a $1$. We denote such a region by $\textup{QH}^{p_1, p_2, \ldots, p_n}_{2n,\ell}$. We note that the bottom row triangles force some rhombus in any tiling, we shall remove such rhombi in the figures that follows.

We are interested in the rhombus tilings of these quartered hexagons. We assign a weight of $x_i$ to each ``left-oriented'' rhombus ($\rightr$) that appears in the $i$-th row of a tiling of $\textup{QH}^{p_1, p_2, \ldots, p_n}_{2n,\ell}$, and let all other rhombi that appear in a tiling to have weight $1$. This weighted region is denoted by 
$$
\textup{QH}_{2n,\ell}^{p_1, p_2, \ldots, p_n}(x_1,x_2, \ldots, x_{2n})
$$
and the generating function of its rhombus tilings by 
\[
M(\textup{QH}_{2n,\ell}^{p_1, p_2, \ldots, p_n}(x_1,x_2, \ldots, x_{2n})).
\]
The result which we will use is the following.
\begin{prop}[Theorem 2.8, \cite{ilse-sp}]\label{prop-ilse}
For a partition $\lambda=(\lambda_1, \lambda_2, \ldots, \lambda_n)$ allowing also zero parts, we have
\begin{equation}\label{sp-ht}
    Sp_{2n}(\lambda; x_1, x_2, \ldots, x_n)=M(\textup{QH}_{2n,\lambda_1}^{\lambda_n+1, \lambda_{n-1}+2, \ldots, \lambda_1+n}(x_1, \bar x_1, x_2, \bar x_2, \ldots, x_n, \bar x_n)).
\end{equation}
\end{prop}

The quantity that we are interested in equation \eqref{vsh-3} is
\begin{equation}\label{eq-1a}
    Sp_{4n+2}(n,n-1,n-1, \ldots, 1,1,0,0; x^2, 1, \ldots, 1)= M(\textup{QH}_{4n+2,n}^{1, 2, 4, 5, \ldots, 3n-2, 3n-1, 3n+1}(x^2, \bar x^2, 1, \ldots, 1))
\end{equation}
For ease in the sequel, we use the following notation
\[
\qh_n(x) = \textup{QH}_{4n+2,n}^{1, 2, 4, 5, \ldots, 3n-2, 3n-1, 3n+1}(x^2, \bar x^2, 1, \ldots, 1) \text{ and }  \qh_n= \qh_n(1).
\]
In general, we call something \textit{unweighted}, if it has weight $1$.

\begin{figure}
    \centering
    \includegraphics[width=0.4\textwidth]{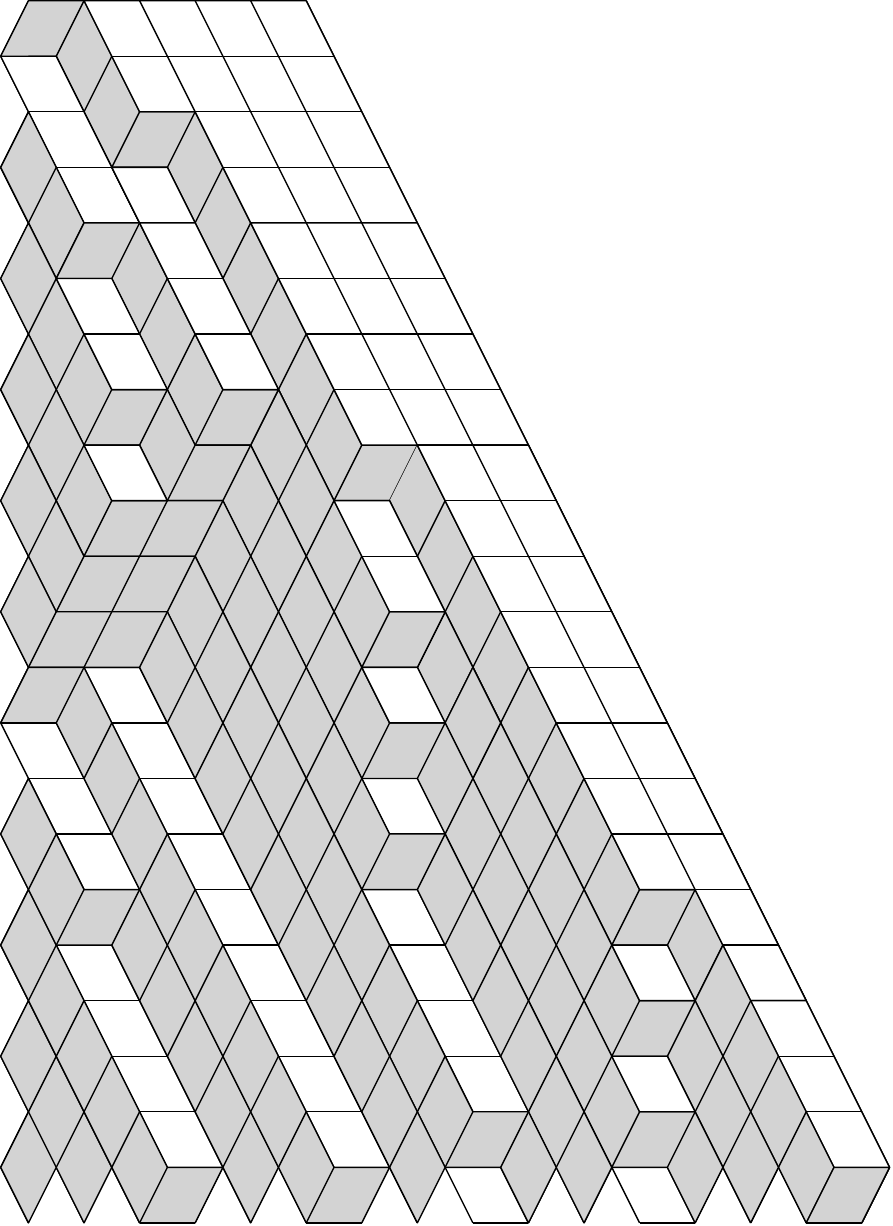}
    \caption{A tiling of $\qh_5$, where $k=i=2$.}
    \label{fig:ht-new-1}
\end{figure}

Note that the weights in our case are manifest only in the first two rows of the quartered hexagon (i.e., one zig-zag portion of the hexagon). There is precisely one vertical rhombus ($\vertr$) in the first row (see Figure~\ref{fig:ht-new-1}), 
say in position $j$ if counted from left and starting with $1$, 
and the rhombi that appear to its right in the first row are all forced to be of the type $\rightr$, while no such rhombus appears left of the vertical rhombus in the first row. Thus, the first row contributes $x^{2(n+1-j)}$ to the weight. In the second row, all rhombi left of the bottom triangle of that fixed vertical rhombus are of type $\rightr$ and they contribute $x^{2(-j+1)}$ to the weight. Now, right of this vertical rhombus there is precisely one vertical rhombus that has its upper triangle in the second row, say in position $i$, and all rhombi right of this are of type $\rightr$, and they contribute $x^{2(-n-1+i)}$ to the weight. The total weight is then 
$x^{2(i-2j+1)}$. 

\begin{figure}
    \centering
    \includegraphics[width=0.4\textwidth]{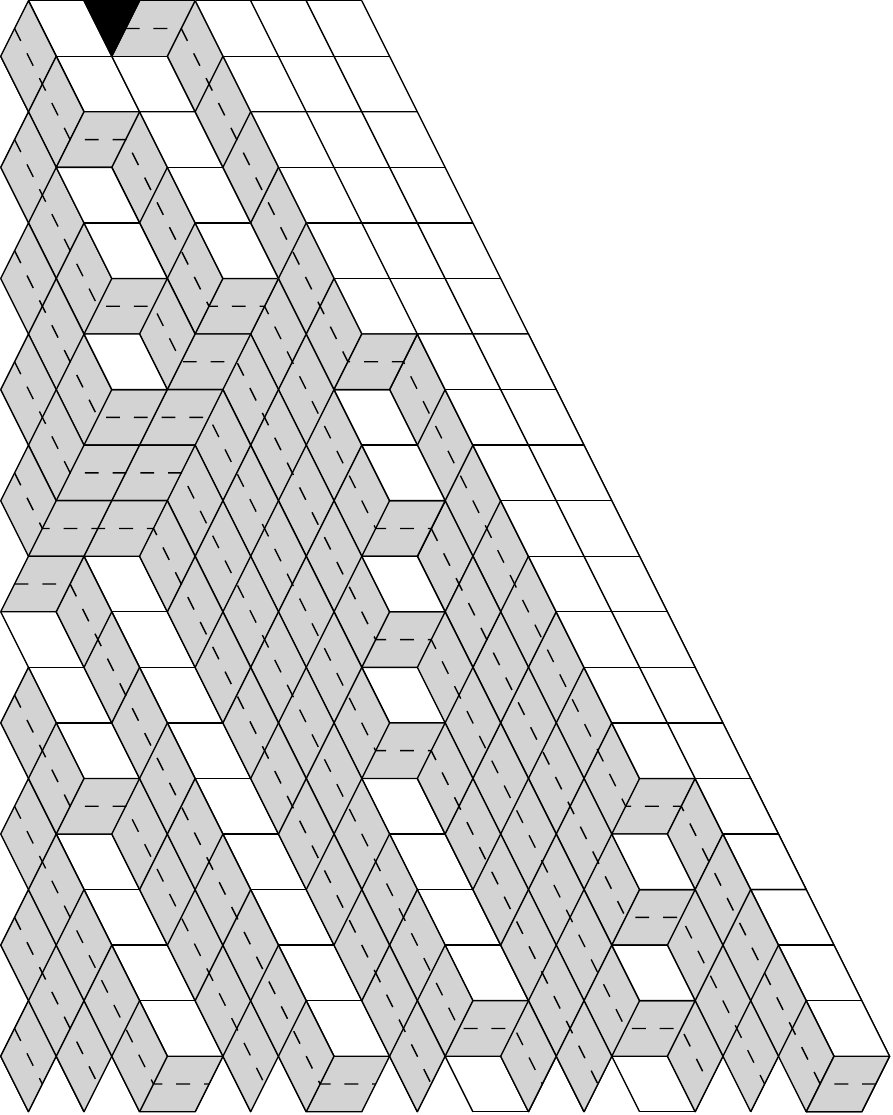}
    \caption{A tiling of $\qh_{5,2}$.}
    \label{fig:ht-new-2}
\end{figure}

We let $Q_{n,i}$ denote the number of tilings of $\qh_n$ with parameter $i$ as described in the previous paragraph and letting $j$ vary in $\{1,2,\ldots,i\}$. Equivalently, letting $\qh_{n,i}$ denote the region obtained from $\qh_n$ by deleting the top two rows  
as well as the $i$-th down-pointing triangle in the (new) top row (see Figure~\ref{fig:ht-new-2} for $\qh_{n,i}$, the dotted lines are to be ignored at the moment), then 
$Q_{n,i}$ is the number of lozenge tilings of $\qh_{n,i}$. 
Then, from the above discussion, we have
\begin{equation}\label{matching}
    M(\qh_n(x))=\sum_{1 \le j \le i \le n+1} Q_{n,i} x^{2i-4j+2}.
\end{equation}

We shall now evaluate $Q_{n,i}$ using the Lindstr\"om-Gessel-Viennot \cite{lind,gessel} technique, where the first step is to transform a rhombus tiling of the region into a family of non-intersecting lattice paths, and then evaluate the total number of such paths by means of a determinant.

We mark the left sides of the up-pointing triangles of the left boundary in our region by $s_1, s_2, \ldots, s_{2n}$ from bottom to top, mark the right side of the deleted (black) down-pointing triangle in the top row by $s_{2n+1}$ and mark the left sides of the deleted up-pointing triangles (by forcing) on the bottom boundary of our region by $e_1, e_2, \ldots, e_{2n+1}$ from left to right. Any rhombus tiling of the region can be transformed into a family of non-intersecting lattice paths, where each path starts at $s_j$ and ends in $e_j$. An example of this in shown in Figure \ref{fig:ht-new-2}, where the lattice paths are given by the dotted lines.

\begin{figure}
   \centering
    \includegraphics[width=0.45\textwidth]{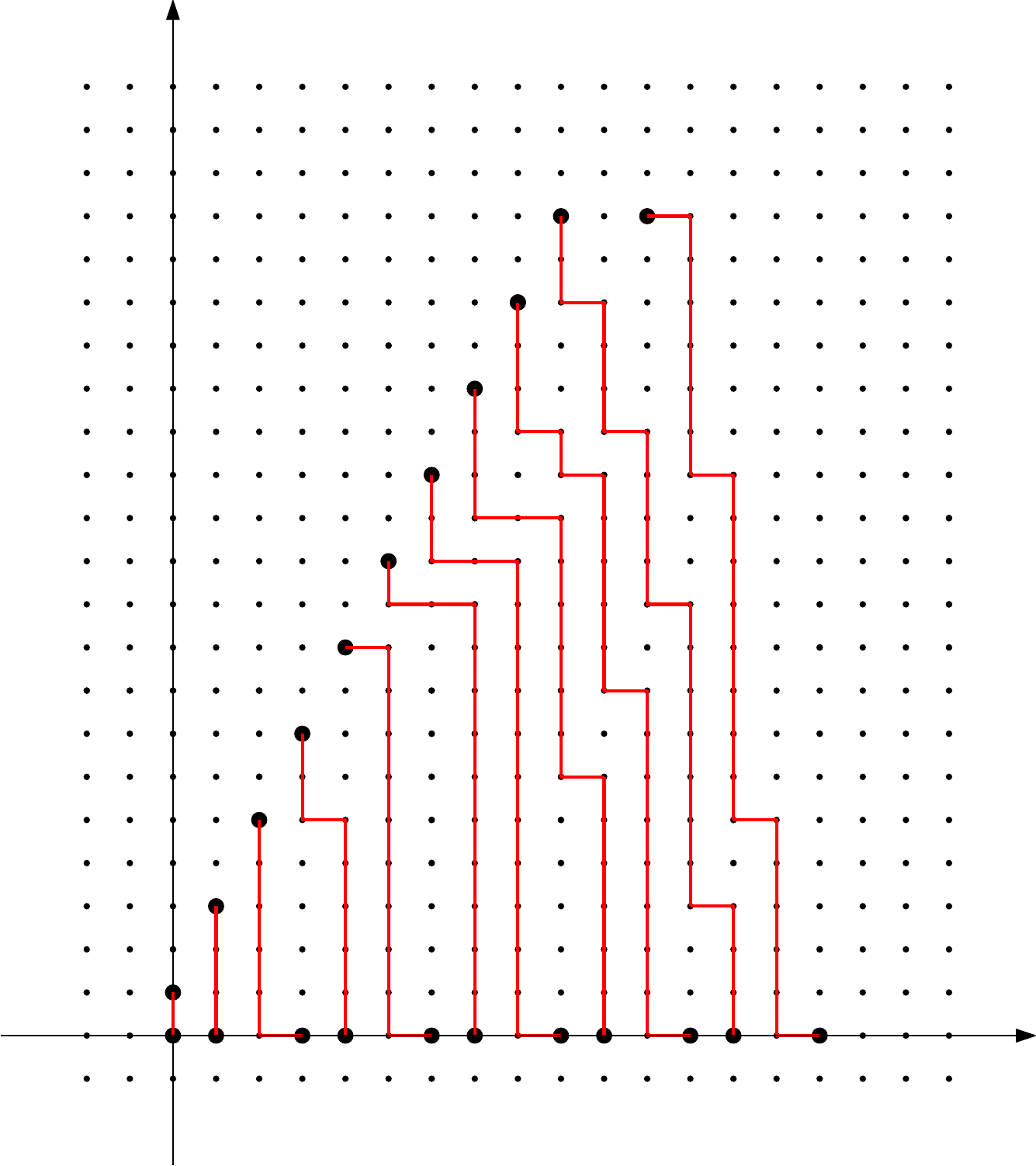}
    \caption{The collection of non-intersecting lattice paths arising from the tiling of $\qh_{5,2}$.}
    \label{fig:ht-4}
\end{figure}

We normalize the obtuse coordinate system in Figure \ref{fig:ht-new-2} and rotate it into the position shown in Figure \ref{fig:ht-4}, where the points $s_j$ have coordinates $(j-1,2j-1)$ for $1\leq j\leq 2n$, $s_{2n+1}$ has coordinates $(2n-1+i,4n-1)$ and $e_j$ have coordinates $(\lfloor (3j-1)/2\rfloor-1,0)$, for $1\leq j\leq 2n+1$. For ease of notation, we denote $\lfloor (3j-1)/2\rfloor$ by $a_j$. Now, we want to find the number of non-intersecting lattice paths beginning at $s_j$ and ending at $e_j$, $1\leq j \leq 2n+1$ with only down and east steps.

By the Lindstr\"om-Gessel-Viennot Lemma \cite{lind,gessel}, the number of such \textit{unweighted} non-intersecting lattice paths is given by the determinant of the $(2n+1) \times (2n+1)$ matrix, whose $(u,v)$-th entry is the number $q_{u,v}$ of lattice paths from $e_u$ to $s_v$. Thus, we have 
\[
q_{u,v}=\begin{cases}
\dbinom{a_u+v-1}{2v-1} & 1\leq u \leq 2n+1, 1\leq v \leq 2n \\
\dbinom{a_u+2n-1-i}{4n-1} & 1\leq u \leq 2n+1, v=2n+1
\end{cases},
\]
and
\[
Q_{n,i}=\det(q_{u,v})_{1\leq u,v\leq 2n+1}.
\]
Expanding the determinant with respect to the $(2n+1)$-th column, we get
\[
Q_{n,i}=\sum_{j=1}^{2n+1}(-1)^{j+1} D_{n,j} \dbinom{a_{j}+2n-1-i}{4n-1},
\]
where
\[
D_{n,j}=\det_{\mycom{1\leq u \leq 2n+1, u \neq j}{1\leq v \leq 2n}} \left(\dbinom{a_u+v-1}{2v-1}\right).
\] 
Taking out common factors and then swapping the $v$-th column with the $(2n-v+1)$-th column for every $v$, we get that $D_{n,j}$ equals
\begin{multline}\label{det-d}
    (-1)^{n(2n-1)}\frac{\prod_{i=1}^{2n+1}a_i}{a_j \prod_{i=1}^{2n}(2i-1)!}  \\ \times
    \det_{\mycom{1 \le i \le 2n+1, i \not= j }{ 1\leq k\leq 2n}} ((a_i-2n+k)(a_i-2n+k+1)\ldots (a_i-1)(a_i+1)\ldots (a_i+2n-j)).
\end{multline}

In order to evaluate this determinant, we need the following lemma, due to Krattenthaler.
\begin{lem}[Lemma 4, \cite{kratt-det}]\label{kratt-lemma}
Let $x_1, x_2, \ldots, x_r$ and  $y_2, y_3, \ldots, y_r$ be indeterminates, and $c$ be a constant. Then we have
\[
\det\left(\prod_{k=j+1}^{r}(x_i-y_k-c)(x_i+y_k)\right)_{1\leq i,j \leq r}=\prod_{1\leq i<j \leq r}(x_j-x_i)(c-x_i-x_j).
\]
\end{lem}

We take $c=0$, $x_i=a_i$ and $y_j=2n-j+1$ in Lemma \ref{kratt-lemma} and evaluate the expression in \eqref{det-d} to get
\[
D_{n,j}= \frac{2\prod\limits_{i=1}^{2n+1}a_i}{\prod\limits_{i=1}^{2n}(2i-1)!} \frac{\prod\limits_{1\leq p <q \leq 2n+1}(a_q-a_p)(a_p+a_q)}{\prod\limits_{q=1}^{2n+1}(a_j+a_q)\prod\limits_{q=j+1}^{2n+1}(a_q-a_j)\prod\limits_{q=1}^{j-1}(a_j-a_q)}.
\]

Therefore, we now have
\begin{equation}
    \label{gf}
        M(\qh_n(x))=\sum_{i=1}^{n+1}\left[\left(\frac{x^{4i}-1}{x^{2i-2}(x^4-1)}\right)\sum_{j=1}^{2n+1}(-1)^{j+1} D_{n,j} \dbinom{a_{j}+2n-1-i}{4n-1} \right].
\end{equation}

\begin{rem}
If we do not delete the first two rows of $\qh_n$, then taking a different family of non-intersecting lattice paths with a similar weighing scheme, rather than the ones used here, we would have arrived at the following generating function:
\begin{multline*}
 M(\qh_n(x))=\det\left(\left(x^2+\frac{1}{x^2}\right)\left\{\binom{4n}{2n-3j+i}-\binom{4n}{2n-3j-i}\right\}+\binom{4n}{2n-3j+i+1}\right.\\
     \left.-\binom{4n}{2n-3j-i-1}+\binom{4n}{2n-3j+i-1}-\binom{4n}{2n-3j-i+1} \right)_{1\leq i,j\leq n}.
\end{multline*}
\end{rem}

\bibliographystyle{alpha}


\begin{thebibliography}{MRR83}

\bibitem[AB19]{AyyerBehrend}
Arvind Ayyer and Roger~E. Behrend.
\newblock Factorization theorems for classical group characters, with
  applications to alternating sign matrices and plane partitions.
\newblock {\em J. Combin. Theory Ser. A}, 165:78--105, 2019.

\bibitem[ABF16]{AyyerBehrendFischer}
Arvind Ayyer, Roger~E. Behrend, and Ilse Fischer.
\newblock Extreme diagonally and antidiagonally symmetric alternating sign
  matrices of odd order.
\newblock {\em arXiv:1611.03823v1 [math.CO]}, 2016.

\bibitem[AD10]{aval-duchon}
Jean-Christophe Aval and Philippe Duchon.
\newblock Enumeration of alternating sign matrices of even size
  (quasi-)invariant under a quarter-turn rotation.
\newblock {\em Electron. J. Combin.}, 17(1):Research Paper 51, 20, 2010.

\bibitem[AF19]{ilse-sp}
Arvind Ayyer and Ilse Fischer.
\newblock Bijective proofs of skew schur polynomial factorizations.
\newblock {\em arXiv:1905.05226v1 [math.CO]}, 2019.

\bibitem[AR13]{AyyerRomik}
Arvind Ayyer and Dan Romik.
\newblock New enumeration formulas for alternating sign matrices and square ice
  partition functions.
\newblock {\em Adv. Math.}, 235:161--186, 2013.

\bibitem[Beh13]{behrend}
Roger~E. Behrend.
\newblock Multiply-refined enumeration of alternating sign matrices.
\newblock {\em Adv. Math.}, 245:439--499, 2013.

\bibitem[BFK17]{bfk}
Roger~E. Behrend, Ilse Fischer, and Matja\v{z} Konvalinka.
\newblock Diagonally and antidiagonally symmetric alternating sign matrices of
  odd order.
\newblock {\em Adv. Math.}, 315:324--365, 2017.

\bibitem[Bre99]{proofs}
David~M. Bressoud.
\newblock {\em Proofs and confirmations}.
\newblock MAA Spectrum. Mathematical Association of America, Washington, DC;
  Cambridge University Press, Cambridge, 1999.
\newblock The story of the alternating sign matrix conjecture.

\bibitem[CLP98]{CohnLarsenPropp}
Henry Cohn, Michael Larsen, and James Propp.
\newblock The shape of a typical boxed plane partition.
\newblock {\em New York J. Math.}, 4:137--165, 1998.

\bibitem[Duc08]{duchon}
Philippe Duchon.
\newblock On the link pattern distribution of quarter-turn symmetric {FPL}
  configurations.
\newblock In {\em 20th {A}nnual {I}nternational {C}onference on {F}ormal
  {P}ower {S}eries and {A}lgebraic {C}ombinatorics ({FPSAC} 2008)}, Discrete
  Math. Theor. Comput. Sci. Proc., AJ, pages 331--342. Assoc. Discrete Math.
  Theor. Comput. Sci., Nancy, 2008.

\bibitem[Fis07]{ilse-1}
Ilse Fischer.
\newblock A new proof of the refined alternating sign matrix theorem.
\newblock {\em J. Combin. Theory Ser. A}, 114(2):253--264, 2007.

\bibitem[Fis09]{FischerOpFormulaVSASM}
Ilse Fischer.
\newblock An operator formula for the number of halved monotone triangles with
  prescribed bottom row.
\newblock {\em J. Combin. Theory Ser. A}, 116(3):515--538, 2009.

\bibitem[Fis11]{ilse-refined}
Ilse Fischer.
\newblock Refined enumerations of alternating sign matrices: monotone
  {$(d,m)$}-trapezoids with prescribed top and bottom row.
\newblock {\em J. Algebraic Combin.}, 33(2):239--257, 2011.

\bibitem[Fis16]{ilse-2}
Ilse Fischer.
\newblock Short proof of the {ASM} theorem avoiding the six-vertex model.
\newblock {\em J. Combin. Theory Ser. A}, 144:139--156, 2016.

\bibitem[FR09]{ilse-romik}
Ilse Fischer and Dan Romik.
\newblock More refined enumerations of alternating sign matrices.
\newblock {\em Adv. Math.}, 222(6):2004--2035, 2009.

\bibitem[GV85]{gessel}
Ira Gessel and G\'{e}rard Viennot.
\newblock Binomial determinants, paths, and hook length formulae.
\newblock {\em Adv. in Math.}, 58(3):300--321, 1985.

\bibitem[KR10]{romik}
Matan Karklinsky and Dan Romik.
\newblock A formula for a doubly refined enumeration of alternating sign
  matrices.
\newblock {\em Adv. in Appl. Math.}, 45(1):28--35, 2010.

\bibitem[Kra99]{kratt-det}
C.~Krattenthaler.
\newblock Advanced determinant calculus.
\newblock {\em S\'{e}m. Lothar. Combin.}, 42:Art. B42q, 67, 1999.
\newblock The Andrews Festschrift (Maratea, 1998).

\bibitem[Kup96]{asm-kuperberg}
Greg Kuperberg.
\newblock Another proof of the alternating-sign matrix conjecture.
\newblock {\em Internat. Math. Res. Notices}, (3):139--150, 1996.

\bibitem[Kup02]{KuperbergRoof}
Greg Kuperberg.
\newblock Symmetry classes of alternating-sign matrices under one roof.
\newblock {\em Ann. of Math. (2)}, 156(3):835--866, 2002.

\bibitem[Lin73]{lind}
Bernt Lindstr\"{o}m.
\newblock On the vector representations of induced matroids.
\newblock {\em Bull. London Math. Soc.}, 5:85--90, 1973.

\bibitem[MRR83]{asm-conj}
W.~H. Mills, David~P. Robbins, and Howard Rumsey, Jr.
\newblock Alternating sign matrices and descending plane partitions.
\newblock {\em J. Combin. Theory Ser. A}, 34(3):340--359, 1983.

\bibitem[Oka06]{OkadaCharacters}
Soichi Okada.
\newblock Enumeration of symmetry classes of alternating sign matrices and
  characters of classical groups.
\newblock {\em J. Algebraic Combin.}, 23(1):43--69, 2006.

\bibitem[Pro94]{proctor}
Robert~A. Proctor.
\newblock Young tableaux, {G}el\cprime fand patterns, and branching rules for
  classical groups.
\newblock {\em J. Algebra}, 164(2):299--360, 1994.

\bibitem[Rob91]{robbins-1}
David~P. Robbins.
\newblock The story of {$1,2,7,42,429,7436,\cdots$}.
\newblock {\em Math. Intelligencer}, 13(2):12--19, 1991.

\bibitem[Rob00]{robbins-2}
David~P. Robbins.
\newblock Symmetry classes of alternating sign matrices.
\newblock {\em arXiv:math/0008045 [math.CO]}, 2000.

\bibitem[RS04]{RazumovStroganov}
A.~V. Razumov and Yu.~G. Stroganov.
\newblock On a detailed enumeration of some symmetry classes of
  alternating-sign matrices.
\newblock {\em Teoret. Mat. Fiz.}, 141(3):323--347, 2004.

\bibitem[RS06a]{rs-refined}
A.~V. Razumov and Yu.~G. Stroganov.
\newblock Bethe roots and refined enumeration of alternating-sign matrices.
\newblock {\em J. Stat. Mech. Theory Exp.}, (7):P07004, 12, 2006.

\bibitem[RS06b]{rs-odd-half}
A.~V. Razumov and Yu.~G. Stroganov.
\newblock Enumeration of odd-order alternating-sign half-turn-symmetric
  matrices.
\newblock {\em Teoret. Mat. Fiz.}, 148(3):357--386, 2006.

\bibitem[RS06c]{rs-odd-quarter}
A.~V. Razumov and Yu.~G. Stroganov.
\newblock Enumeration of odd-order alternating-sign quarter-turn symmetric
  matrices.
\newblock {\em Teoret. Mat. Fiz.}, 149(3):395--408, 2006.

\bibitem[Sta86]{stanley}
Richard~P. Stanley.
\newblock A baker's dozen of conjectures concerning plane partitions.
\newblock In {\em Combinatoire \'enum\'erative ({M}ontreal, {Q}ue.,
  1985/{Q}uebec, {Q}ue., 1985)}, volume 1234 of {\em Lecture Notes in Math.},
  pages 285--293. Springer, Berlin, 1986.

\bibitem[Str06]{stroganov-izergin}
Yu.~G. Stroganov.
\newblock The {I}zergin-{K}orepin determinant at a cube root of unity.
\newblock {\em Teoret. Mat. Fiz.}, 146(1):65--76, 2006.

\bibitem[Tsu98]{Uturn}
Osamu Tsuchiya.
\newblock Determinant formula for the six-vertex model with reflecting end.
\newblock {\em J. Math. Phys.}, 39(11):5946--5951, 1998.

\bibitem[Zei96a]{doron}
Doron Zeilberger.
\newblock Proof of the alternating sign matrix conjecture.
\newblock {\em Electron. J. Combin.}, 3(2):Research Paper 13, approx.\ 84,
  1996.
\newblock The Foata Festschrift.

\bibitem[Zei96b]{zeil}
Doron Zeilberger.
\newblock Proof of the refined alternating sign matrix conjecture.
\newblock {\em New York J. Math.}, 2:59--68, electronic, 1996.

\end{thebibliography}

\end{document}